\documentclass[a4paper,reqno,10pt,oneside]{amsart}
\usepackage{amsmath,geometry}
\usepackage{graphicx}
\usepackage{mathrsfs}
\usepackage{amssymb,amsmath, amsfonts, amsthm}
\usepackage{latexsym}
\usepackage{color}
\usepackage[dvips]{epsfig}
\usepackage{graphicx}

\allowdisplaybreaks[1]

\numberwithin{equation}{section}

\renewcommand{\O}{\operatorname{O}}

\renewcommand{\(}{\left(}
\renewcommand{\)}{\right)}
\renewcommand{\[}{\left[}
\renewcommand{\]}{\right]}

\newtheorem{theorem}{Theorem}[section]
\newtheorem{proposition}[theorem]{Proposition}
\newtheorem{lemma}[theorem]{Lemma}

\newtheorem{remark}[theorem]{Remark}

\renewcommand{\le}{\leqslant}
\renewcommand{\ge}{\geqslant}

\renewcommand{\d }{\delta }

\newcommand{\g }{\gamma }

\renewcommand{\l }{\lambda}

\renewcommand{\t}{\theta}
\renewcommand{\O}{\Omega}

\newcommand{\U}{\mathcal{U}}
\newcommand{\m}{\lambda}

\newcommand{\pp}{\partial}

\newcommand{\beq}{\begin{equation}}
\newcommand{\eeq}{\end{equation}}
\newcommand{\beqs}{\begin{equation*}}
\newcommand{\eeqs}{\end{equation*}}
\newcommand{\beqn}{\begin{eqnarray}}
\newcommand{\eeqn}{\end{eqnarray}}
\newcommand{\beqns}{\begin{eqnarray*}}
\newcommand{\eeqns}{\end{eqnarray*}}
\newcommand{\bdoc}{\begin{document}}
\newcommand{\edoc}{\end{document}}
\newcommand{\be}{\begin{enumerate}}
\newcommand{\ee}{\end{enumerate}}
\newcommand{\bdescr}{\begin{description}}
\newcommand{\edescr}{\end{description}}
\newcommand{\ba}{\begin{array}}
\newcommand{\ea}{\end{array}}
\newcommand{\intR}{\int_{\mathbb R^N}}
\newcommand{\R}{\mathbb R^N}
\newcommand{\e}{\varepsilon}
 \renewcommand{\(}{\left(}
 \newcommand{\la}{\lambda}
\renewcommand{\)}{\right)}

\newcommand{\ve}{ \varepsilon}
\newcommand{\rr}{ \mathbb{R}}

\renewcommand{\[}{\left[}
\renewcommand{\]}{\right]}
\newenvironment{Proof}{\noindent{\bf Proof.}}{\hfill$\Box$\\[2mm]}

\begin{document}

\title[]{Large mass boundary condensation patterns in the stationary Keller-Segel system}

\author{Manuel del Pino}
\address[Manuel del Pino]{Departamento de Ingenieria Matem\'atica
Facultad de Ciencias Fisicas y Matem\'aticas,Universidad de Chile
Casilla 170, Correo 3, Santiago, Chile}
\email{delpino@dim.uchile.cl}

\author{Angela Pistoia}
\address[Angela Pistoia] {Dipartimento SBAI, Sapienza Universit\`{a} di Roma, via Antonio Scarpa 16, 00161 Roma, Italy}
\email{angela.pistoia@uniroma1.it}

\author{Giusi Vaira}
\address[Giusi Vaira] {Dipartimento SBAI, Sapienza Universit\`{a} di Roma, via Antonio Scarpa 16, 00161 Roma, Italy}
\email{vaira.giusi@gmail.com}


\begin{abstract}
We consider the boundary value problem
\begin{equation*}
\left\{
\begin{array}{lr}
-\Delta u + u =\lambda e^u, \quad \mbox{in}\,\, \Omega\\
\partial_{\nu}u =0\qquad\qquad \mbox{on}\,\, \partial \Omega
\end{array}
\right.
\end{equation*}
where $\Omega$ is a bounded smooth domain in $\mathbb R^2$,
$\lambda>0$ and $\nu$ is the inner normal derivative at $\partial\Omega$. This problem
is equivalent to the stationary Keller-Segel
 system from chemotaxis.\\
We establish the existence of a solution $u_\la$  which exhibits a sharp boundary layer along the entire  boundary $\partial\Omega$ as $\lambda\to 0$. These solutions have large mass in the sense that $ \int_\Omega \la e^{u_\la}  \sim |\log\la|.$

\end{abstract}

\maketitle

 \section{Introduction and statement of the main result}
Chemotaxis is one of the simplest mechanisms for aggregation of biological species. The term
refers to a situation where organisms, for instance bacteria, move towards high concentrations
of a chemical which they secrete. A basic model in chemotaxis was introduced by Keller
and Segel \cite{KS}. They considered an advection-diffusion system consisting of two coupled
parabolic equations for the concentration of the considered species and that of the chemical
released, represented, respectively, by positive quantities $v(x, t)$ and $u(x, t)$ defined on a
bounded, smooth domain $\Omega$ in $\mathbb R^N$ under no-flux boundary conditions. The system reads
\beq\label{pb1}
\left\{
\ba{lr}
\frac{\partial v}{\partial t}=\Delta v -\nabla \cdot (v\nabla u)\qquad &\mbox{in}\quad\O\\
\tau\frac{\partial u}{\partial t}=\Delta u -u +v\,\,\qquad &\mbox{in}\quad \O\\
%u, v>0 \,\,\qquad &\mbox{in}\quad \O\\
\frac{\partial u} {\partial\nu }= \frac{\partial v} {\partial\nu } =0 \qquad\quad   &\qquad\mbox{on}\quad\partial\O,
\ea
\right.
\eeq
where $\nu$ denotes the unit inner normal to $\partial\Omega$.
Steady states of \eqref{pb1} are the positive solutions of the system
\beq\label{pb3}
\left\{
\ba{lr}
\Delta v -\nabla\cdot (v\nabla u)=0\qquad &  \mbox{in}\quad\O\\
\Delta u -u+v=0\qquad\quad   & \mbox{in}\quad\O\\
\frac{\partial u}{\partial\nu}=  \frac{\partial v}{\partial\nu}  =0 \qquad\quad\   &\qquad\mbox{on}\quad \partial\O.
\ea
\right.
\eeq
Problem \eqref{pb3} can be reduced to a  scalar equation. Indeed, testing the first equation against $(\ln v -u)$,
an integration by parts shows that a solution of \eqref{pb3}  satisfies the relation
$$
\int_\Omega  v|\nabla (\ln v  -  u) |^2 = 0
$$
and hence $v= \l e^u$ for some positive constant $\l$, and thus $u$ satisfies the equation

\beq\label{pb}
\left\{
\ba{lr}
-\Delta u + u =\l e^u, \qquad \mbox{in}\quad \O,\\
\ \quad \frac{\partial u}{\partial\nu}  =0\qquad\ \ \qquad\mbox{on}\quad\partial\O.
\ea
\right.
\eeq
Reciprocally, a solution to problem \eqref{pb} produces one of \eqref{pb3} after setting $v= \l e^u$. In this paper we consider problem \eqref{pb} when
$\O\subset\mathbb R^2$ is a bounded domain with smooth boundary and $\l>0$ is a small parameter. By integrating both sides of the equation we see that a necessary condition for existence is
$\l<1$.

The analysis of problems \eqref{pb1}, \eqref{pb3} and their corresponding versions in entire space $\mathbb R^ 2$, has a long history, starting with the work by Childress and Percus \cite{cp}.
 The analysis of the steady state problem \eqref{pb} for small $\l$  started with Schaaf \cite{schaaf} in the one-dimensional case. Existence of a radial solution when $\Omega$ is a ball, generating a spike shape at the origin when $\l\to 0$ was established by Biler \cite{biler}.  The shape of an unbounded family of solutions $u_\l$ with uniformly bounded masses
$$
\limsup_{\l\to 0^+ } \int_\O \l e^{u_\l} < +\infty $$ was established in \cite{suzuki1,wei}. As in the classical analysis by Brezis and Merle \cite{bm},  blow-up of the family is found to occur at most   on a finite number of points $\xi_1,\ldots,\xi_k\in \Omega$, $\xi_{k+1},\ldots,  \xi_{k+l}\in \partial\Omega$. More precisely, in the sense of measures,
\begin{equation}\label{jjj}
-\Delta u_\l + u_\l = \l e^{u_\l}  \rightharpoonup   \sum_{i=1}^ k 8\pi \delta_{\xi_i}  + \sum_{i=k+1}^ {k+l} 4\pi \delta_{\xi_{i}}
\end{equation}
as $\l\to 0$. Here $\delta_\xi$ denotes the Dirac mass at the point $\xi$. Correspondingly,
away from those points the leading behavior of $u_\m$ is given by
\begin{equation}\label{kk}
u_\l(x) \to  \sum_{i=1}^ k 8\pi G(x,\xi_i)  + \sum_{i=k+1}^ {k+l} 4\pi G(x, \xi_{i})
\end{equation}
where $G(\cdot,\xi)$ is the Green function for the problem
\beq
\left\{
\ba{lr}
-\Delta G + G =  \delta_\xi, \qquad \mbox{in}\quad \O,\\
\ \quad \frac{\partial G}{\partial\nu}=0\qquad\ \ \qquad\mbox{on}\quad\partial\O.
\ea
\right.
\eeq

For each given non-negative numbers $k$ and $l$, a solution  $u_\l$ with the properties \eqref{jjj} and \eqref{kk} for suitable points $\xi_i$ is proven to exist in \cite{DW}.
Near each point $\xi= \xi_i$  the leading concentration behavior is given by  $$u_\l(x)\sim \omega(|x-\xi|)$$ where $\omega$ is a radially symmetric solution of the equation
\begin{equation}\label{kkk}-\Delta \omega = \l e^\omega  \quad\hbox{in }\mathbb R^2,\end{equation}
namely a function of the form
$$\omega(r) = \ln \frac{8\delta^ 2} { (\delta^ 2 + r^2 )^ 2} - \ln \lambda . $$
where $\delta$ is a suitable scalar dependent on $\m$ and the point  $\xi$.

Since $u_\l$ is uniformly bounded away from the
 points $\xi_i$, this forces for the parameter $\delta$ to satisfy  $\delta^2 \sim \lambda$. We observe that
 all solutions $\omega$ of \eqref{kkk} satisfy $$\int_{\mathbb R^2} \l e^\omega = 8\pi.$$ Thus, consistently with \eqref{jjj}, masses are  quantized as
 \beq \label{11}
 \int_\Omega \l e^{u_{\l}} \to 4\pi (2k + l)  . \eeq

\medskip
A natural question is that of  analyzing  of solutions with large mass, namely solutions $u_\l$ of \eqref{pb3}  with   $$\int_\Omega \l e^{u_{\l}} \to  +\infty \quad\hbox{as }\l \to 0.$$
It is natural to seek for solutions with  property which concentrate not just at points but  on a {\em larger-dimensional set}.
The purpose of this paper is to prove the existence of a family of solutions to \eqref{pb3}
with a {\em boundary condensation property},  exhibiting a boundary layer behavior along the entire
$\partial\Omega$. These solutions satisfy
$$
\lim_{\l \to 0 }  \frac 1 {| \ln \l | } \int_\Omega \lambda e^{u_\l}  >0.  $$
Let us formally derive the asymptotic shape of these solutions.
Let us parametrize points of space in a sufficiently small neighborhood of $\partial\Omega$ in the form
$
x= \gamma(\theta) + y\nu(\theta),
$
where $\gamma(\theta) $ is a parametrization by $\theta$, arclength of $\partial\Omega$, and $\nu(\theta)$ a corresponding unit inner normal, so that $\dot \nu(\theta) = -\kappa(\theta) \dot \gamma(\theta)$, where $\kappa$ designates inner normal curvature.   We get the following expansion for the Euclidean Laplacian in these coordinates
$$
\Delta   =  \partial_{yy} +  \frac 1 { 1-\kappa(\theta)y} \frac{\partial}{\partial{\theta}} \left ( \frac 1 { 1-\kappa(\theta)y} \frac {\partial } {\partial \theta} \right ) - \frac {\kappa(\theta)}  { 1-\kappa(\theta)y}\frac { \partial } {\partial y}
$$
The solution we look for has a boundary layer, thus large derivatives along the normal and a comparatively smooth behavior along the tangent direction. It is then
reasonable to take near $\partial\O$ as a first approximation of a solution $u(\theta,y)$ of the equation \eqref{pb}
 a solution of the ordinary differential equation
\begin{equation}\label{w0}
w^{''}_\mu  + \lambda e^{w_\mu} = 0, \quad w_\mu '(0) = 0,
\end{equation}
which is
\begin{equation} w_\mu (y)    {-\ln\lambda}= w(y/\mu) - 2\ln \mu    -  \ln \lambda \label{dd},\ \end{equation}
where $w(y) = \ln 4\frac {e^{\sqrt2 y}} {\(1+e^{\sqrt2 y}\)^2} $ and the concentration parameter   $\mu $ satisfies
 $$\mu(\theta) = \e \hat \mu   {_\e}(\theta)\sim \e \hat \mu_0(\theta).$$ Here $\e=\e(\l)$ is a small positive number which we shall choose below
and $\hat \mu_0(\theta)$ is a uniformly positive and bounded  smooth function.

Let $\varphi \in C(\bar \Omega)$ compactly supported near the boundary of $\Omega.$ A direct computation yields

$$
    \e  \int_{\Omega} \l e^ {w_\mu} \varphi   =  \sqrt 2  \int_{\partial\Omega}  \varphi \hat \mu_0^{-1} \, d\theta  + O(\e)\ ,
$$
since $\int\limits_0^{+\infty}e^{w(y)}dy=\sqrt2.$
Thus,
$$
\e \l e^ {w_\mu} \rightharpoonup   \sqrt 2\hat\mu_0^{-1}\delta_{\partial\Omega}
$$
where $\delta_{\partial\Omega}$ is the Dirac measure on the curve ${\partial\Omega}$.

Then we expect that, globally,
$ \sqrt 2\ \mathcal U = \e u_\lambda $ satisfies approximately
$$ -\Delta  \mathcal U +  \mathcal U =     \hat \mu_0^{-1} \delta_{\partial\Omega} ,
$$
which means in the limit
$$- \Delta  \mathcal U +  \mathcal U  =0 \hbox{ in } \Omega, \quad  \partial_\nu  \mathcal U=  - \hat \mu_0^{-1} \hbox{ on } \partial\Omega.
$$
Now, from our ansatz \eqref{dd}, we should have that close to the boundary
 $$ \sqrt 2\ \mathcal U(\theta,y)\approx  \e w(y/\mu) - 2\e \ln \e\hat\mu    - \e \ln \lambda  $$
and hence, in particular
$$
 \sqrt 2\ \mathcal U (\theta,0) \approx- \e \ln \lambda -2\e\ln\e
$$
By maximum principle and $ \partial_\nu  \mathcal U=  - \hat \mu^{-1}_0<0$
the latter relation is consistent in the limit if the constant $\e \ln \lambda$ approaches a    {negative} number. If we choose $\mathcal U  =1$ on the boundary of $\O, $ then we take $\e$  such that
$$-\e  \ln \lambda -2\e \ln\e  \approx  \sqrt2$$ so that
$$
\e  \approx -\frac {\sqrt2 }{\ln \la}   .
$$

Hence the limiting $ \mathcal U$ equals $ \mathcal U_0$, the unique solution of the problem
 \begin{equation}\label{G}
- \Delta  \mathcal U_0 +  \mathcal U_0  =0 \hbox{ in } \Omega, \quad    \mathcal U_0 = 1 \hbox{ on } \partial\Omega.
\end{equation}
We observe that by maximum principle and Hopf's Lemma, we have that  $\partial_\nu  \mathcal U_0 <0$, and hence
this fixes our choice of $\hat\mu_0(\t)$ as
$$
\hat\mu_0(\t) = -\frac 1 {\partial_\nu  \mathcal U_0}  \quad\hbox{on } \partial \Omega.
$$

Our main result asserts the existence of a solution with exactly the profile above for all $\la$ sufficiently small which remains suitably away from a sequence of critical small values where certain resonance phenomenon occurs.

\begin{theorem}\label{principale}
   {Suppose that $\Omega$ is a smooth bounded domain of $\mathbb R^2$. Then there exists a sequence of positive small numbers $\lambda=\lambda_m$ converging to $0$ as $m\to+\infty$ such that the problem \eqref{pb} has a solution $u_\l$}
%$$\mathfrak I_m:=[\mathfrak a_m,\mathfrak b_m],\ m\in\mathbb N,\ \hbox{with}\ 0<\mathfrak a_m <\mathfrak b_m\to0\ \hbox{as}\    {m\to+\infty}$$
%such that if $\lambda\in(0,\lambda_0)$ and
%$\lambda\not\in\cup_{m\in\mathbb N} \mathfrak I_m$
 %there exists a solution $u_\la$ to problem \eqref{pb}
such that
$$
0< \lim_{\m\rightarrow 0} \frac 1 {|\ln \la| } \int_\O \lambda e^{u_\m}\, dx <+\infty.
$$
Moreover, if $\e_\lambda=\e_{\l_m}$ is the parameter defined by
\begin{equation}\label{el}
\ln\frac4{\e_\lambda^2}-\ln\lambda={\sqrt 2\over\e_\lambda}
\end{equation}
then
$$
\lim_{\m\rightarrow 0}\e_\m u_\m = \sqrt 2\ \mathcal U_0 \qquad C^0-\mbox{uniformly on compact sets of }\,\O
$$
and, in the sense of measures,
$$
\e_\la  \m e^ {u_\la } \rightharpoonup -\sqrt2\  \partial_\nu  \mathcal U_0\ \delta_{\partial\O}.
$$

\end{theorem}

We actually believe that problem \eqref{pb} has a solution which concentrates along the entire boundary, also in the higher-dimensional case  $\Omega \subset \mathbb R^N$ with $N\ge3.$ This fact has been established in the radial case, when $\Omega$ is a ball, in \cite{pv}.
\\\\

   { Remark \ref{4.3} below assures the existence of small numbers $\l>0$ for which the problem \eqref{pb} has a solution with the desidered behavior. In fact, a more general condition on $\e_\l$ (and then on $\l$) defined as in \eqref{el} is provided there. This type of condition, known as non-resonance condition, were imposed to establish the presence of higher dimensional concentration patterns without rotational symmetries in several works in the literature,
 starting with the pioneering works by Malchiodi and Montenegro \cite{malchiodi-montenegro-cpam,malchiodi-montenegro-duke}, who prove existence
 of a concentrating solution $u_\e$  along the boundary for the classical Neumann problem
 \beq\label{ccc}\e^2 \Delta u - u + u^p=0 \quad\hbox{in } \Omega, \quad \pp_\nu u = 0 \hbox{ on } \pp\Omega \eeq
 with $p>1$. See also  \cite{delpino-kowalczyk-wei-cpam},
\cite{mahmoudi-malchiodi-adv}, \cite{malchiodi-gafa} for related results.
}

A major difference between our problem and \eqref{ccc} is that the limiting profile is highly localized in the sense that the limiting solution has an exponentially sharp boundary layer $O( e^{-\frac d \e})$ where $d$ designates distance to the boundary. Instead, in our setting the interaction with the inner part of the domain is much stronger.  %which is reflected in the asymptotic inner profile  $\mathcal{U}_0 $. 
The interaction inner-outer problem makes the improvement of approximations considerably more delicate. The construction of an inverse for the approximate linearized operator is in fact quite
different because of the presence of slow decay elements in the kernel of the asymptotic linearization.\\\\

The proof of our result  relies on  an  infinite-dimensional form of Lyapunov-Schmidt reduction.\\
We look for a solution to \eqref{pb}  of the form   $U_\lambda+\Phi_\lambda$ where $U_\lambda$, {\em the main term}, is a suitably constructed first approximation and 
$\Phi_\lambda$ is {\em the remainder term}.
Then Problem \eqref{pb} can be rewritten as 
\begin{equation}\label{pro1}
   {+}L(\Phi_\lambda)=S _\lambda+N (\Phi_\lambda)\ \hbox{in}\ \Omega,
\end{equation}
where
\begin{equation}\label{L}
L(\Phi ):=\Delta \Phi -\Phi+\lambda e^{U_\lambda}\Phi,
\end{equation}
\begin{equation}\label{S}
S_\lambda(U_\l):=-\Delta U_\lambda +U_\lambda-\lambda e^{U_\lambda}
\end{equation}
and
\begin{equation}\label{N}
N(\Phi ):=    {-}\lambda e^{U_\lambda}\left[e^\Phi-1-\Phi\right].
\end{equation}
The strategy consists of finding an accurate first approximation $U_\lambda$ (Section \ref{1}) so that the error term $S_\lambda(U_\l)$ be small in a suitably chosen norm (Section \ref{2}). Then 
an invertibility theory for associated linearized operator $L$ (Section \ref{4}) allows to solve equation \eqref{pro1} for term $\Phi_\lambda $ via a fixed point argument (Section \ref{3}).
\\
The main term $U_\lambda$ looks  like $w_\mu-\ln\lambda$  close to the boundary, with $w_\mu$ defined in \eqref{dd}
  solves the ODE \eqref{w0}
and concentration parameter $\mu:=\mu(\lambda)$ approaches 0 as   $\lambda$ goes to 0.
The profile of $U_\lambda$ in the inner part of the domain looks like $\tau\ \mathcal U_0$ where $ \mathcal U_0$ solves the Dirichlet boundary problem
\eqref{G}
and the dilation parameter $\tau:=\tau(\lambda)$ approaches $+\infty$ as $\lambda$ goes to 0.
The concentration parameter $\mu(\lambda)$ and the dilation parameter $\tau(\lambda)$ have to be chosen so that the two profiles match accurately close to the boundary. This is the most delicate part of the paper and it is carried out in sub-section \ref{10}.

 \section{The main term}\label{1}

\subsection{The problem close to the boundary}
  Let us parametrize $\partial\O$ by the arc length  $$\gamma(\theta):=(\gamma_1(\theta), \gamma_2(\theta)),\qquad \theta\in [0,\ell]$$ where $\ell:=|\partial\O|$. The tangent vector and the inner normal vector to the point $\gamma(\theta)\in\partial\O$ are given by $$\tau(\theta):=(\dot{\gamma}_1(\theta), \dot{\gamma}_2(\theta));\qquad \nu(\theta):=(-\dot{\gamma}_2(\theta), \dot{\gamma}_1(\theta))$$ respectively.\\
If $\delta>0$ is small enough,   let $$\mathcal D_{\delta}:=\left\{x\in \O\,\,\,:\,\,\, {\rm dist}(x, \partial\O)\leq \delta\right\}$$ be a neighbourhood of the curve $\partial\O$.\\ Then for any $x\in \mathcal D_{\delta}$ there exists a unique $(\theta, y)\in [0, \ell]\times [-\delta, 0]$ such that $$ x=\gamma(\theta)+ y \nu(\theta)=(\gamma_1(\theta)-y\dot\g_2(\t), \g_2(\t)+y\dot\g_1(\t)).$$
We remark that in these coordinates the points of the boundary take the form $(\t, 0)$.
If $  u(\t, y)$ is a function defined in $[0, \ell]\times [-    \d, 0]$ we can define the function
 $u(x)=  u (\t( x), y( x))$ (we use the same symbol for sake of simplicity) for $x\in\mathcal D_{\delta}$  and hence close to the boundary the equation   \eqref{pb} takes the form
\beq\label{pbcoo}
\left\{\begin{aligned}
&-\displaystyle\frac{1}{(1-y\kappa(\t))^2}\partial_{\t\t}^2  u-\partial^2_{yy} u -\displaystyle\frac{y\dot\kappa(\t)}{(1-y\kappa(\t))^3}\partial_\t  u +\displaystyle\frac{   {\kappa(\t)}}{1-y\kappa(\t)}\partial_y  u + u =\lambda e^{ u}\quad \mbox{in}\,\,\, \mathcal D_{\delta},\\
\\
&\partial_ {   {y}} u(\t,0)=0
\end{aligned}\right.
\eeq
where $\kappa(\t)$ is the curvature at the point $\g(\t)\in\partial\O$. \\

It is useful to introduce the spaces $C^0_\ell(\mathbb R)$ and $C^2_\ell(\mathbb R)$   of $\ell-$periodic $C^0-$functions and $C^2-$functions, respectively.

\subsection{The  scaled problem close to the boundary }\label{preliminari}

Now, let us introduce an extra parameter $\varepsilon:=\varepsilon_\lambda$  such that
\beq\label{relelambda}
\ln\frac{4}{\e_\lambda^2}-\ln\lambda=\frac{\sqrt{2}}{\e_\lambda},\qquad \hbox{i.e.}\quad \lambda={4\over\e_\m^2}e^{-{\sqrt2\over\e_\m}}.
\eeq
It is easy to check that $\e_\lambda\rightarrow 0 $ as $\lambda\rightarrow 0$. We agree that in the following we will use indifferently the two parameters $\e$ and $\l$ to get the necessary estimates. \\
 Moreover, let us choose  the concentration parameter $\mu(\theta):=\mu(\lambda,\theta)$ in \eqref{dd}  as
\beq\label{mudef}
\mu(\t):=\e \hat\mu (\t), \hbox{where}\ \hat\mu (\t):=\hat\mu_\e(\t)\in C^2_\ell(\mathbb R).
\eeq
The function $\hat\mu  $ will be defined in Lemma \ref{ue}.
\\
Finally, let us set
\beq\label{tildemu0}
\hat\mu_0(\t):=-\frac{1}{\partial_\nu \U_0}\Big|_{\partial\O}=-\frac{1}{\partial_y \U_0(\t, 0)}.
\eeq
We note that by maximum principle and Hopf's lemma, $ \mu_0 $ is a strictly positive $C^2-$function.\\

Now, let us scale problem \eqref{pbcoo}.
In $\mathcal D_{\delta}$ it is natural to consider the change of variables
\beq\label{cambio}
 u(\t, y)=\tilde u\left(\frac{\t}{\e}, \frac{y}{\mu}\right),\ \hbox{with}\ \tilde u=\tilde u(s, t).
\eeq
It is clear that
$$(\t, y)\in \mathcal D_{\delta}\ \hbox{if and only if}\ (s,t)\in \left[0, \frac{\ell}{\e}\right]\times\left[-\frac{\delta}{\mu}, 0\right].$$

   {Let $\tilde u=\tilde u(s, t)$, then we can   compute
\begin{equation*}
\begin{aligned}
\partial_\t   u &=\e^{-1}\partial_s\tilde u -\dot\mu\mu^{-1} t \partial_t  \tilde u \\
\partial_{\t\t}^2  u &=\e^{-2} \partial_{ss}^2  \tilde u -2\e^{-1}\dot\mu \mu^{-1}t\partial_{st}^2\tilde u-\ddot\mu\mu^{-1} t  \partial_{t}\tilde u +2\dot\mu^2\mu^{-2} t \partial_t \tilde u +\dot\mu^2\mu^{-2}t^2 \partial_{tt}^2\tilde u\\
\partial_y u &= \mu^{-1} \partial_t \tilde u\\
\partial_{yy}^2 u &= \mu^{-2}\partial^2_{tt}\tilde u
\end{aligned}
\end{equation*}}
where the dot stands for the derivative with respect to $\t$.\\
Hence, problem   \eqref{pbcoo}   can be written as
\beq\label{pbcoocambio}
\left\{
\begin{array}{lr}
\hat\mu_0^2(\e s) \partial_{ss}^2 \tilde u +\partial_{tt}^2 \tilde u +\tilde{\mathcal A}(\tilde u)+\lambda\mu^2 e^{\tilde u}=0\qquad \mbox{in}\,\,\, \mathcal C_{\delta},\\\\
\partial_{   {t}} \tilde u =0 \qquad\qquad\qquad\qquad\qquad \mbox{on}\,\,\, \partial \mathcal C_\delta \cap \{t=0\}\\\\
\tilde u\left(s+\displaystyle\frac{\ell}{\e}, t\right)= \tilde u(s, t),\quad \hbox{i.e. $\tilde u$ is $\frac\ell\e-$periodic in $s$}\\\\
\end{array}
\right.
\eeq
where $\mathcal C_\delta:=\mathbb R^+\times \left[-\frac{\delta}{\mu}, 0\right]$ and
the linear operator $\tilde {\mathcal A}$ is defined by
\begin{eqnarray}\label{tildeA}
\nonumber
\tilde{\mathcal A}(\tilde u)&:=& \underbrace{\left[\frac{\hat\mu^2}{\left(1-\mu t \kappa(\e s)\right)^2}-\hat\mu_0^2\right]}_{b_0(s, t)}\partial_{ss}^2 \tilde u +   {\underbrace{\frac{\dot\mu^2 t^2}{\left(1-\mu t \kappa(\e s)\right)^2}}_{b_1(s, t)}}\partial_{tt}^2 \tilde u \\
\nonumber
&&   {- \underbrace{\frac{2\e^{-1}\mu\dot\mu t}{\left(1-\mu t \kappa(\e s)\right)^2}}_{b_2(s, t)}}\partial^2_{st} \tilde u 
+\underbrace{\frac{\mu^3\e^{-1}t \dot\kappa(\e s)}{\left(1-\mu t \kappa(\e s)\right)^3}}_{b_4(s, t)}\partial_s\tilde u-\mu^2 \tilde u\\
\nonumber
&&+   {\underbrace{\left[-\frac{\ddot\mu \mu t}{\left(1-\mu t \kappa(\e s)\right)^2}-\frac{\mu^2\dot\mu t^2 \dot\kappa(\e s)}{\left(1-\mu t \kappa(\e s)\right)^3}+\frac{2\dot\mu^2 t}{(1-\mu t \kappa(\e s))^2}-\frac{\mu\kappa(\e s)}{1-\mu t \kappa(\e s)}\right]}_{b_3(s, t)}}\partial_t \tilde u\\.
\end{eqnarray}
It is important to point out that the linear operator $\tilde{\mathcal A}$ is a perturbation term since all $b_i$'s are uniformly small when $\lambda$ is small (because of \eqref{mudef}).

\subsection{A linear theory close to the boundary}
Let us read the first order term   of  $u_\lambda$ close to the boundary in the scaled variables: since $u_\lambda$ looks like $w_\mu-\ln\lambda$   where the one-dimensional bubble  $w_\mu$ is defined in  \eqref{dd},  it turns out that the first order term of $\tilde u_\lambda$ is nothing but   $w-\ln\lambda$ where $w\equiv w_1,$ namely
\beq\label{w}
w(t):=\ln  4\frac{e^{\sqrt 2 t}}{\(1+e^{\sqrt 2 t}\)^2}, \quad t\in\mathbb R
\eeq
which solves
\beq\label{pbw}
  w^{''}+e^w=0, \qquad \mbox{in}\quad \mathbb R.
\eeq

Therefore, it is important to develop a linear theory for the linear operator $\mathcal L$ which comes from the linearization of equation \eqref{pbcoocambio} around the bubble $w-\ln\lambda,$
namely
\begin{equation}\label{operatoreL}
\mathcal L (\tilde\phi ):=\hat\mu_0^2 \partial_{ss}^2 \tilde\phi  +\partial^2_{tt}\tilde\phi +\tilde{\mathcal A}(\tilde\phi )+ e^w \tilde\phi.
\end{equation}

In order to study $\mathcal L $,  an important role is played by the linear operator \begin{equation}\label{oppallino}\hat   {\mathcal L}  (\tilde\phi ):=\partial_{tt}^2 \tilde\phi  + e^w \tilde\phi \end{equation}
which  is nothing but the linearized operator around $w$ of equation \eqref{pbw}.\\
 \begin{lemma}\label{zetas}   Let us consider
 the  associated linearized eigenvalue problem
$$  \hat  {\mathcal L}(\tilde\phi )=\Lambda\tilde\phi \quad \mbox{in}\,\, \mathbb R.$$
 \begin{itemize}
 \item[(i)] $\Lambda=0$ is an eigenvalue with associated eigenfunctions    $Z_1(t)=2+t   w'(t)$ and $ Z_2(t)= w'(t)=\sqrt 2 \frac{1-e^{\sqrt 2 t}}{1+e^{\sqrt 2 t}}.$  We point out that  $Z_1$ behaves like a constant at infinity and that $Z_2$ is not a bounded function.
\item[(ii)] There exists  a positive eigenvalue $\Lambda_1$ with associated radial, positive and bounded eigenfunction $Z_0=Z_0(t)$ with $L^2-$norm equal to one. Moreover, $Z_0$ decays exponentially at infinity as $O\left(e^{-{\sqrt\Lambda_1 |t|}}\right).$
\end{itemize}
\end{lemma}
\begin{Proof}(i) has been proved in \cite{g3}. (ii) can be proved arguing as in Section 3 in \cite{dmp}.
\end{Proof}
We consider the following projected problem: {\em given    a bounded function $h$, which is ${\ell\over\varepsilon}-$periodic in $s, $ find s bounded ${\ell\over\varepsilon}-$periodic function $c_0(s)$ and $\tilde\phi$ such that}
\beq\label{sistemaL}
\left\{
\begin{array}{lr}
\mathcal L(\tilde\phi )= h + c_0(s) Z_0(t)\qquad \mbox{in}\,\, \mathcal C_{\delta},\\\\
\partial_{   {t}} \tilde\phi =0\qquad \mbox{on}\,\,\, \partial\mathcal C_{\delta}   {\cap \{ t=0\}}, \\\\
\tilde\phi \left(s+\displaystyle\frac{\ell}{\e}, t\right)=\tilde\phi(s, t)\\\\
\displaystyle\int_{-\frac{2\d}{\mu}}^0 \tilde\phi (s, t) Z_0(t)\, dt =0\,\,\, \qquad \forall\,\, s\in \mathbb R^+ .
\end{array}
\right.
\eeq

In Section \ref{4}, we will establish existence and a priori estimates for problem \eqref{sistemaL} in the following norms:
\begin{equation}\label{normas}
\|\phi \|_{*}:=\sup_{\mathcal C_{\delta}}(1+|t|^\sigma)|\phi|+\sup_{\mathcal C_{\delta}}(1+|t|^{\sigma+1})|\nabla\phi |,\  \|h\|_{**}:=\sup_{\mathcal C_{\delta}}(1+|t|^{\sigma+2})|h|,\ \hbox{for}\
\sigma\in (0, 1).\end{equation} More precisely, we prove that

\begin{proposition}\label{inv}
There exist $\lambda_0>0$ and a constant $C>0$,  such that for any $\lambda\in(0,\lambda_0)$ and for any $h$ with $\|h\|_{**}<+\infty$, there exists a unique $\phi=\mathcal T(h)$ bounded solution of the problem \eqref{sistemaL} such that
\begin{equation}\label{controllonorme}
\|\phi\|_{*}\leq C\|h\|_{**}.
\end{equation}
\end{proposition}

\subsection{The main term close to the boundary}

The function $w_\mu-\ln\l$ is the main term of the approximated solution close to the boundary. We need to add   some correction terms, which improve the main term.\\ More precisely, we let
\begin{equation}\label{acb}u_\lambda(\t, y)=\underbrace{w_\mu(\t, y)-\ln\l + \alpha_\mu(\t, y)}_{\hbox{$1^{st}-$order}}+\underbrace{v_\mu(\t, y)+\beta_\mu(\t, y)}_{\hbox{$2^{nd}-$order}}+\underbrace{z_\mu(\t, y)}_{\hbox{$3^{rd}-$order}}+ \underbrace{e_0^\e(\t) Z_0^\mu(y)}_{\hbox{unknown!}}\end{equation} where
\begin{itemize}
\item $\alpha_\mu(\t, y)$ is defined in Lemma \ref{lemmaalpha},
\item $v_{\mu}(\t, y)$ is defined in Lemma \ref{lemmav} and $\beta_\mu(\t, y)$ is defined in Lemma \ref{lemmabeta},
\item $z_\mu(\t, y)$ is defined in Lemma \ref{lemmaz}
\item $Z_0^\mu(y)=Z_0(\frac{y}{\mu})$, where $Z_0$ is defined in Lemma \ref{zetas} and the function $e_0^\e(\t)$  is defined as follows
 \beq\label{eje}
 e_0^\e(\t)=   {\e^{\frac 32}} e_0(\t)\ \hbox{with}\ e_0\in C^2_\ell(\mathbb R) .
 \eeq

 We point out that the function $e_0$ is unknown: it is playing the role of one parameter and it will be chosen in  Section \ref{fap} as  solution of an ordinary differential equation. We assume that $e_0$ has uniformly bounded  $\|\cdot\|_\e-$norm, i.e.
 \beq\label{normae}
 \|e_0\|_\e:=\|\e^2\ddot{e_0}\|_\infty+\|\e \dot e_0\|_\infty +\|e_0\|_\infty\le M_0,
 \eeq
 for some large fixed number $M_0$.
 \end{itemize}

 The first term we have to add is a sort of projection of the function $w_\mu$, namely the function $\alpha_\mu$ given in the next lemma.\\
\begin{lemma}\label{lemmaalpha}
\begin{itemize}
\item[(i)] The Cauchy problem
\beq\label{alpha}
\left\{
\ba{lr}
-\partial_{yy}^2\alpha_\mu+\displaystyle\frac{\kappa(\t)}{1-y\kappa(\t)}\partial_y \alpha_\mu=
-\frac{\kappa(\t)}{1-y\kappa(\t)}\partial_y w_\mu -w_\mu+\ln\l    {+\frac{1}{(1-y\kappa(\t))^2}\partial^2_{\t\t} w_\mu}\\\\
\alpha_\mu(\t, 0)=\partial_y \alpha_\mu(\t, 0)=0
\ea
\right.
\eeq
has the solution
\begin{equation*}
\begin{aligned}
\alpha_\mu(\t, y)&=-\int_0^y \frac{1}{1-\sigma \kappa(\t)}\int_0^{\sigma}(1-\rho\kappa(\t))\[-\frac{\kappa(\t)}{1-\rho\kappa(\t)}\partial_y w_\mu( \rho)-w_\mu(\rho)+\ln\l\]\,d\sigma\,d\rho\\
&   {-\int_0^y \frac{1}{1-\sigma\kappa(\t)}\int_0^\sigma \frac{1}{1-\rho\kappa(\t)}\partial^2_{\t\t} w_\mu (\rho)\, d\rho\, d\sigma}
\end{aligned}
\end{equation*}
\item[(ii)] For any $(\t, y)\in \mathcal D_{2\delta}\setminus \mathcal D_\delta$ it holds:

\begin{eqnarray*}
\alpha_{\mu}(\t, y)&:=&(\kappa(\t)\ln 4) y+\frac{y^2}{2}\left[\kappa^2(\t)\ln 4   {+\frac{d^2}{d\t^2}\ln\hat\mu^2 }-\frac{\sqrt 2}{\mu}\kappa(\t) +\(\ln\frac{4}{\mu^2}-\ln\l\)\right]+\\
&&+\frac{y^3}{6}\left[-\frac{2\sqrt 2}{\mu}\kappa^2(\t)-\frac{\sqrt 2}{\mu}+\kappa(\t)\left(\ln\frac{4}{\mu^2}-\ln\l\right)   {+\frac{\sqrt 2}{\e}\frac{d^2}{d\t^2}\frac{1}{\hat\mu}}\right]\\
&&+O\( |y|^{   {3}}\)+O\(  \frac{|y|^4}{\e}\).
\end{eqnarray*}
\item[(iii)] Moreover, via the change of variables $\t=\e s$ and $y=\mu t,$ the function  $\tilde\alpha_\mu(s, t):=\alpha_\mu (\e s, \mu t)$ solves the problem
    \beq\label{tildealpha}
\left\{\begin{aligned}
&-\partial_{tt}^2\tilde\alpha_\mu+\displaystyle\frac{\mu(\e s) \kappa(\e s)}{1- t \mu(\e s)   \kappa(\e s)}\partial_t \tilde\alpha_\mu=
-\frac{\mu(\e s)\kappa(\e s)}{1-t\mu(\e s) \kappa(\e s)}\partial_t w - \mu^2\(\ln\frac{1}{\mu^2}-\ln\lambda+w \) \\
&   {+\frac{1}{(1-\mu t\kappa(\t))^2}\left[\hat\mu^2 \partial_{ss}^2 \ln \frac{1}{4\hat\mu^2}+t\partial_t w \left(-\frac{\ddot\mu}{\mu}+\frac{2\dot\mu^2}{\mu^2}\right)\mu^2+\dot\mu^2 t^2\partial^2_{tt} w\right]}\\\\
&\tilde\alpha_\mu(\e s, 0)=\partial_t \tilde\alpha_\mu(\e s, 0)=0.\\
\end{aligned}
\right.
\eeq
\item[(iv)] The following expansion holds
$$\tilde{\alpha}_\mu(s, t):=\alpha_\mu(\e s, \mu t)=\mu \alpha_1(\e s, t)+\mu^2 \alpha_2(\e s, t)+O( \e^3  t^{   {4}}),\qquad    {\mbox{for}\,\, |t| \le \frac{2\delta}{\mu} }$$
where
\beq\label{a1}
\alpha_1(\e s, t)=\kappa(\e s)\int_0^t w(\sigma)\, d\sigma+\frac{\sqrt 2}{2}\hat\mu (\e s) t^2\eeq
and
\beq\label{a2}
\begin{aligned}
\alpha_2(\e s,t)&= \kappa^2 (\e s)\int_0^t \sigma w(\sigma)\,d\sigma+   {\frac{t^2}{2}\ln\frac{1}{\hat\mu^2(\t)}} \\
&+\int_0^t\int_0^\sigma \(w(\rho)-\ln 4\)\,d\rho\,d\sigma +   {\frac{\sqrt 2}{6}\hat\mu \kappa(\t) t^3+\frac{d^2}{d\t^2}\ln\hat\mu^2 \frac{t^2}{2}}.\end{aligned}\eeq

\end{itemize}
\end{lemma}
\begin{Proof}
We argue as in s Lemma 3.1 of \cite{pv}.
\end{Proof}
\\
 Now, let us construct the  second order term of our approximated solution.  \begin{lemma}\label{lemmav}
\begin{itemize}\item[(i)]
There exists $v$ solution of the linear problem ($\alpha_1$ is given in \eqref{a1})
\beq\label{vuffa}-\partial_{yy}^2v-e^w v =e^w \alpha_1(\t, y)\eeq such that $$v(\t, y)=\nu_1(\t) y +\nu_2(\t)+O(e^{-|y|})\qquad |y|\rightarrow +\infty$$ where
\beq\label{nu1}
\nu_1(\t):=2\kappa(\t)(1-\ln2)+\ln 4 \hat\mu(\t)\eeq
and \beq\label{nu2}
\nu_2(\t):=-\int_{-\infty}^0 \left(\frac{2}{1-e^{\sqrt{2}y}}+\frac{y}
{\sqrt{2}}\right)\alpha_1(\t, y)\partial_y w(y) e^{w(y)}\, dy.\eeq
\item[(ii)]In particular, the function $v_\mu(\theta,y):=\mu v\(\theta,\frac y\mu\)$ solves the problem
$$ -\partial_{yy}^2v_\mu-e^{w_\mu} v_\mu =\mu e^{w_\mu} \alpha_1\(\theta,\frac y\mu\).$$
\item[(iii)]
Moreover, via the change of variables $\t=\e s$ and $y=\mu t,$ the function $\tilde v_\mu(s, t):= v_\mu(\e s,\mu t)=\mu (\e s) v(\e s, t)$  solves the problem
$$-\partial_{tt}^2 \tilde v_\mu - e^w \tilde v_\mu =\mu e^w \alpha_1(\e s, t)  $$
and   the following expansion holds (see \eqref{mudef})
    $$\tilde v_\mu(s, t):= \e \nu_1(\e s) \hat\mu (\e s)t +\e\nu_2(\e s) \hat\mu(\e s) + O\(\e e^{-|t|}\)\ \hbox{as}\    {|t|\to+\infty}.$$  \end{itemize}
\end{lemma}
\begin{Proof}
We apply Lemma \ref{massimo}.
\end{Proof}

As we have done for the function $w_\mu$, we have to add the projection of the function $v_\mu$, namely the function $\beta_\mu$ given in the next lemma.

\begin{lemma}\label{lemmabeta}
\begin{itemize}
\item[(i)] The Cauchy problem ($v_\mu$ is given in Lemma \ref{lemmav})
\begin{equation}\label{cauchybeta}
\left\{
\begin{array}{lr}
-\partial_{yy}^2\beta_\mu+\displaystyle\frac{\kappa(\t)}{1-y\kappa(\t)}\partial_y\beta_\mu=-\displaystyle\frac{\kappa(\t)}{1-y\kappa(\t)}\partial_y v_\mu(\t,y)\\\\
\beta_\mu(\t, 0)=\partial_y \beta_\mu(\t, 0)=0
\end{array}
\right.
\end{equation}
has the solution
$$\beta_\mu(\t, y)=\int_0^y \frac{\kappa(\t)}{1-\sigma\kappa(\t)}\int_0^{\sigma}   {\partial_y v_\mu(\t, \rho)}\, d\rho\,d\sigma.$$
\item[(iv)] For any $(\theta, y)\in \mathcal D_{2\delta}\setminus \mathcal D_\delta$ we have:
$$\beta_\mu(\t, y)= \nu_1(\t) \kappa(\t)\frac{y^2}{2}+O( |y|^3).$$
\item[(iii)] Moreover, via the change of variables $\t=\e s$ and $y=\mu t,$ the function $\tilde\beta_\mu(s, t):=\beta_\mu(\e s, \mu t)$ solves the problem
\begin{equation}\label{cauchybetatilde}
\left\{
\begin{array}{lr}
-\partial_{tt}^2\tilde\beta_\mu+\displaystyle\frac{\mu (\e s)\kappa(\e s)}{1-t\mu(\e s)\kappa(\e s)}\partial_t\tilde\beta_\mu=-\displaystyle\frac{\mu(\e s)\kappa(\e s)}{1-t\mu(\e s)\kappa(\e s)}\partial_t \tilde v_\mu\\\\
\tilde\beta_\mu(\e s, 0)=\partial_t \tilde\beta_\mu(\e s, 0)=0
\end{array}
\right.
\end{equation}
\item[(iv)]
The following expansion holds:
$$\tilde\beta_\mu(s, t):=\beta_\mu(\e s, \mu t)=\mu^2 \beta_1(\e s, t)+O(\e^3  t^3)$$ where
\begin{equation}\label{b1}
\beta_1(\e s, t)=\kappa(\e s)\int_0^t \int_0^{\sigma}    {\partial_y v(\e s, \rho)}\, d\rho\,d\sigma.\eeq

\end{itemize}
\end{lemma}
\begin{Proof}
We argue as in Lemma 3.4 of \cite{pv}.
\end{Proof}
\\
 Finally, we build   the  third order term of our approximated solution.
\begin{lemma}\label{lemmaz}
\begin{itemize}
There exists $z$ solution of the linear problem ($\alpha_1,$   $\alpha_2$ and $\beta_1$ are given in \eqref{a1}, \eqref{a2} and \eqref{b1} respectively, and $v$ is given in Lemma \ref{lemmav}
)
$$-\partial_{yy}^2 z-e^w z =e^w \left[\alpha_2(\t, y)+\beta_1(\t, y)+\frac 12 \left(\alpha_1(\t, y)+ v(\t, y)\right)^2\right]$$ such that $$z(\t, y)=\zeta_1(\t) y +\zeta_2(\t)+O(e^{-|y|})\qquad |y|\rightarrow +\infty$$ where $$\zeta_1(\t):=\frac{1}{\sqrt 2}\int_{-\infty}^0 h(\t, y)\partial_y w(y) e^{w(y)}\, dy$$ and $$\zeta_2(\t):=-\int_{-\infty}^0 \left(\frac{2}{1-e^{\sqrt{2}y}}+\frac{y}
{\sqrt{2}}\right)h(\t, y)\partial_y w(y) e^{w(y)}\, dy,$$
with $$h(\t, y)=\alpha_2(\t, y)+\beta_1(\t, y)+\frac 12 \left(\alpha_1(\t, y)+ v(\t, y)\right)^2.$$
\item[(ii)]In particular, the function $z_\mu(\theta,y):=\mu^2 z\(\theta,\frac y\mu\)$ solves the problem
$$ -\partial_{yy}^2z_\mu-e^{w_\mu} z_\mu =\mu^2 e^{w_\mu} \left[\alpha_2\(\theta,\frac y\mu\)+\beta_1\(\theta,\frac y\mu\)+\frac 12 \left(\alpha_1\(\theta,\frac y\mu\)+ v\(\theta,\frac y\mu\)\right)^2\right].$$
\item[(iii)] Moreover, via the change of variables $\t=\e s$ and $y=\mu t,$ the function $\tilde z_\mu(s, t):=z_\mu(\e s, \mu t)=\mu^2 (\e s) z(\e s, t)$ solves the problem
    $$-\partial_{tt}^2 \tilde z_\mu - e^w \tilde z_\mu =\mu^2 e^w \left[\alpha_2(\e s, t)+\beta_1(\e s, t)+\frac12\left(\alpha_1(\e s, t)+ v(\e s, t)\right)^2\right] $$
   and   the following expansion holds (see \eqref{mudef})
$$\tilde z_\mu(s, t):= \e^2\zeta_1(\e s) \hat\mu^2(\e s) t +\e^2\zeta_2(\e s) \hat\mu^2(\e s) + O\(\e^2 e^{-|t|}\)\ \hbox{as}\    {|t|\to+\infty}.$$
 \end{itemize}

\end{lemma}

\begin{Proof}
We apply Lemma \ref{massimo}.
\end{Proof}

\begin{lemma}\label{massimo}(Lemma 4.1, \cite{g3})
Let $\mathfrak h\in C^0(\mathbb R) $ such that $\int\limits_{\mathbb R}\mathfrak h(y)w'(y)e^{w(y)}dy<+\infty.$ Then the function
$$\mathfrak U(y):=w'(y)\int\limits_0^y{1\over\(w'(\sigma)\)^2}\int\limits_\sigma^ 0\mathfrak h(\tau) w'(\tau)e^{w(\tau)}d\tau d\sigma$$
solves the ordinary differential equation
$$-\mathfrak U''-e^w\mathfrak U=e^w \mathfrak h\ \hbox{in}\ \mathbb R.$$
In particular,
$$\mathfrak U(y)=\mathfrak a y +\mathfrak b +O\(e^{-c|y|}\)\ \hbox{in $C^1(\mathbb R)$ as $y\to-\infty $} $$
where
$$\mathfrak a :={1\over\sqrt2}\int\limits^0_{-\infty}\mathfrak h(\tau) w'(\tau)e^{w(\tau)}d\tau\ \hbox{and}\ \mathfrak b :=-\int\limits^0_{-\infty}\({2\over1-e^{\sqrt2 \tau}}+{\tau\over\sqrt 2}\)\mathfrak h(\tau) w'(\tau)e^{w(\tau)}d\tau$$
and
$$\mathfrak U(y)=\mathfrak c y +\mathfrak d +O\(e^{-c|y|}\)\ \hbox{in $C^1(\mathbb R)$ as $y\to+\infty $} $$
where
$$\mathfrak c :={1\over\sqrt2}\int\limits_0^{+\infty}\mathfrak h(\tau) w'(\tau)e^{w(\tau)}d\tau\ \hbox{and}\ \mathfrak d :=-\int\limits_0^{+\infty}\({2\over1-e^{\sqrt2 \tau}}+{\tau\over\sqrt 2}\)\mathfrak h(\tau) w'(\tau)e^{w(\tau)}d\tau.$$
\end{lemma}

\subsection{How to match the main term close to the boundary with the main term in inner part}\label{10}
The solution  $u_\lambda$ in the inner part of the domain looks like $\tau \mathcal U_0$ where $ \mathcal U_0$ solves \eqref{G}
and the dilation parameter $\tau:=\tau(\lambda)$ approaches $+\infty$ as $\lambda$ goes to 0.
The  function $u_\lambda$ (and its derivative)  built in \eqref{acb}   in a neighborhood of the boundary   has to match  with the   function $\tau(\lambda)\mathcal U_0$
(and its derivative). To this aim it  is necessary  to  choose the dilation parameter $\tau={\sqrt2\over\varepsilon}$ and most of all it is essential to modify the profile of the solution in the inner part of the domain by building a new function $\mathcal U_\varepsilon$  which approaches $\mathcal U_0$ as $\varepsilon$ goes to zero  and such that
its value on the boundary together with the value of its normal derivative coincide with the value of $u_\lambda$ and its normal derivative. The main tool here is the Dirichlet-to-Neumann map and the key ingredient is the choice of the concentration parameter $\hat\mu $ as showed in the next crucial lemma.

\begin{lemma}\label{ue}
There exists $\e_0$ such that for any $\e\in (0, \e_0)$  there exist a function $\hat\mu_\e \in C^2(\partial\O)$ and a solution   $\mathcal U_\e$ to the problem
 \beq\label{pbfuoriep1}
\left\{
\ba{lr}
-\Delta \U_\e+\U_\e=0, \qquad\mbox{in}\quad\O,\\
\U_\e= 1- {\varepsilon\over\sqrt2}\(\ln\hat\mu_\e^2   {-}\varepsilon\hat\mu_\e\nu_2   {-\e^2 \hat \mu^2 \zeta_2(\t)} \) \qquad\mbox{on}\quad \partial\O,\\
\partial_\nu\U_\e=-\frac{1}{\hat\mu_\e}+\frac{\e}{\sqrt 2}\(2\kappa +\hat\mu_\e\ln 4   {+\e \hat\mu \zeta_1(\t)}\)\qquad\mbox{on}\quad \partial\O.
\ea
\right.
\eeq
 Moreover ($\hat\mu_0$ is given in   \eqref{tildemu0})
    \beq\label{convtildemu}
 \hat\mu_\e = \hat\mu_0 +O(\e) \quad \hbox{in $C^1(\partial\O)$ as $\e\rightarrow 0$}
 \eeq
and
  \beq\label{ueo}\mathcal U_\e =\mathcal U_0+O\(\e\) \quad \hbox{in $C^2(\overline\O)$  as $\e\rightarrow 0$}.\eeq
\end{lemma}
\begin{Proof}
Let us apply the Dirichlet-to-Neumann map, which maps the value on $\partial \O$  of a  harmonic function $\mathcal U$ to the value of its normal derivative $\partial_\nu \mathcal U$ on $\partial \O,$ i.e. $F\(\mathcal U_{|_{\partial\O}}\)=\partial_\nu\mathcal U.$
Therefore, we are going to find a function $\hat\mu\in C^2(\partial\O)$ such that
\beq\label{dton}
   F\( 1- {\varepsilon\over\sqrt2}\(\ln\hat\mu^2   {-}\varepsilon\hat\mu\nu_2    {-\e^2 \hat \mu^2 \zeta_2(\t)}\) \)=-\frac{1}{\hat\mu}+\frac{\e}{\sqrt 2}\(2\kappa +\hat\mu\ln 4   {+\e\hat\mu\zeta_1(\t)}\).
   \eeq
Let
       $$H(\e, \hat\mu)=F\( 1- {\varepsilon\over\sqrt2}\(\ln\hat\mu^2   {-}\varepsilon\hat\mu\nu_2    {-\e^2 \hat \mu^2 \zeta_2(\t)}\)\)+\frac{1}{\hat\mu}-\frac{\e}{\sqrt 2}\(2\kappa(\t)+\hat\mu\ln4+   {\e \hat\mu \zeta_1(\t)}\).$$
       We have that
   $$H(0, \hat\mu_0)= F(1)+\frac{1}{\hat\mu_0}=0,\ \hbox{ since $ F(1)=\partial_\nu \U_0$ and $\hat\mu_0=-{1\over\partial_\nu \U_0}$ (see  \eqref{tildemu0}).}$$ Moreover $$\frac{\partial H}{\partial \hat\mu}(0, \hat\mu_0)=-\frac{1}{\hat\mu_0^2}\neq 0.$$ Hence by the Implicit Function Theorem, there exists a unique $\hat\mu=\hat\mu_\e(\t)\in C^2(\partial \O)$ such that $H(\e, \hat\mu_\e)=0$, namely \eqref{dton} holds.
  Estimates \eqref{convtildemu} and \eqref{ueo}   follow by elliptic standard regularity theory.

 \end{Proof}

\begin{lemma}
Let $\mathcal U_\varepsilon $ be given in Lemma \ref{ue}.
Then there exists $\e_0$ such that for any $\e\in(0,\e_0)$
\beq\label{u1u3}
    { u_{ \l}(\t, y)-{\sqrt2\over\varepsilon}\U_\e(\t, y)=  O\( \e |y|^2\)+O\( \frac{|y|^4}{\e}\)\ \hbox{uniformly in}\ \mathcal D_{2\delta}\setminus \mathcal D_\delta}
 \eeq
 and
 \beq\label{u1primomenou3primo}
   { \partial_y \[u_{ \l}(\t, y)-{\sqrt2\over\varepsilon}  \U_\e(\t, y)\]= O(\e |y|)+O\( \frac{|y|^3}{\e}\)\ \hbox{uniformly in}\ \mathcal D_{2\delta}\setminus \mathcal D_\delta.}
 \eeq
 \end{lemma}

\begin{Proof}
 Let us prove the estimate \eqref{u1u3}. The proof of \eqref{u1primomenou3primo} is similar. \\
 Let  $ \mathcal U$ be a generic harmonic function, namely
 \begin{equation}\label{armo}-\Delta \U +\U =0  \ \mbox{in}\ \O.
 \end{equation}
 Then the expansion of $\frac{\sqrt2}\varepsilon\mathcal U$ on the boundary reads as
  \begin{equation}\label{glu1}
\frac{\sqrt2}\varepsilon\mathcal U(\t, y)=\frac{\sqrt{2}}{\e}\[\U (\t, 0)+y\partial_y \U (\t, 0)+\frac{y^2}{2}\partial_{yy}^2\U (\t, 0)+\frac{y^3}6\partial_{yyy}^3\U (\t, 0) \] +O\(\frac{|y|^4}{\e}\).\end{equation}
Now, let us write  the   expansion of the function $u_\l$
  close to the boundary.
In
  $\mathcal D_{2\delta}\setminus \mathcal D_\delta$ we get $$w_\mu (y)-\ln\l= \ln \frac{4}{\mu^2}-\ln\l-\frac{\sqrt{2} y}{\mu}+
 O(e^{-\sqrt{2}\frac{|y|}{\mu}})={\sqrt2\over\varepsilon} -\ln\hat\mu^2-\frac{\sqrt{2} y}{\mu}+
 O(e^{-\sqrt{2}\frac{|y|}{\mu}}),$$ because $\mu=\varepsilon\hat\mu$ and \eqref{relelambda} holds. Therefore,
 by  Lemmas \ref{lemmaalpha}, \ref{lemmav}, \ref{lemmabeta} and \ref{lemmaz} we deduce
  \begin{eqnarray}\label{glu2}
 u_{  \l}(\t, y) &=&w_\mu(y) -\ln\l+\alpha_\mu(\t, y)+v_\mu(\t, y)+\beta_\mu(\t, y)+z_\mu(\t, y)+ e_0^\varepsilon (\t)Z_0^\mu(y)\nonumber\\
&=&{\sqrt2\over\varepsilon}\underbrace{\[1- {\varepsilon\over\sqrt2}\(\ln\hat\mu^2   {-}\varepsilon\hat\mu\nu_2(\t)   {-\e^2 \hat \mu^2 \zeta_2(\t)}\)\]}_{\sim\ \U (\t, 0)}\nonumber\\
&&+\frac{\sqrt 2}{\e} y \underbrace{\[-\frac{1}{\hat\mu}+\frac{\e}{\sqrt 2}(2\kappa(\t)+\hat\mu\ln 4    {+\e \hat\mu \zeta_1(\t)})\]}_{\sim\ \partial_y\U (\t, 0)}\nonumber\\
&&+\frac{\sqrt 2}{\e} \frac{y^2}{2}\underbrace{\[1-\frac{\kappa(\t)}{\hat\mu}   {+\frac{\e}{\sqrt 2}\left(\frac{d^2}{d\t^2}\ln \hat\mu^2-\ln\hat\mu^2+2\kappa^2(\t)+\kappa(\t)\hat\mu \ln 4\right)}\]}_{\sim\ \partial_{yy}^2\U (\t, 0)}\nonumber\\
&&+\frac{\sqrt 2}{\e} \frac{y^3}{6}\underbrace{\[-\frac1{\hat\mu}-\frac{2}{\hat\mu}\kappa^2(\t)+\kappa(\t)   {+\frac{d^2}{d\t^2}\frac{1}{\hat\mu}}\]}_{\sim\ \partial_{yyy}^3\U (\t, 0)}\nonumber\\ & &+ O\(   { |y|^3}\)+O\(\frac{|y|^4}{\e}\)+O(e^{-{c}\frac{|y|}{\varepsilon}}),
 \end{eqnarray}
 for some $c>0.$
Let us compare \eqref{glu1} with \eqref{glu2}: the first four terms have to be equal!
In particular, it means that we have to find an harmonic function $\mathcal U$ such that the value of  $\mathcal U$ and the value of its normal derivative $\partial _\nu \mathcal U$ on the boundary have to be equal to
$\[1- {\varepsilon\over\sqrt2}\(\ln\hat\mu^2   {-}\varepsilon\hat\mu\nu_2    {-\e^2 \hat \mu^2 \zeta_2(\t)}\)\]$ and $\[-\frac{1}{\hat\mu}+\frac{\e}{\sqrt 2}(2\kappa+\hat\mu\ln 4   {+\e\hat\mu \zeta_1(\t)})\]$, respectively.
This is done in Lemma \ref{ue}.
Therefore, let us  replace  in \eqref{glu1} and  \eqref{glu2} the generic armonic function $\mathcal U$ with the function  $\mathcal U_\e$ which solves problem \eqref{pbfuoriep1}. The first two terms coincide. Now, let us check what happens with the higher order terms, namely terms which  involve the second and third derivatives of $\mathcal U_\e$.
The function $\mathcal U_\e$ solves equation \eqref{pbfuoriep1} which in a neighborhood of the boundary reads as
 \beq\label{ue1}-{1\over (1-y\kappa)^2}\partial^2_{\theta\theta}\mathcal U_\e-\partial^2_{yy}\mathcal U_\e-{y\dot \kappa\over (1-y\kappa)^3}\partial _{\theta }\mathcal U_\e+
 {  \kappa\over (1-y\kappa) }\partial _{y }\mathcal U_\e+\mathcal U_\e=0\ \hbox{in}\ \mathcal D_{2\delta}. \eeq
%On the other hand,  on the boundary $y=0$ and so the boundary condition  satisfy
 %\beq\label{ue2}\mathcal U_\e(\theta,0)=1+O(\varepsilon),\   \partial_y\mathcal U_\e(\theta,0)=-{1\over\hat\mu_0(\t)}+O(\varepsilon)\ \hbox{and}\ \partial^2_{\theta\theta}\mathcal U_\e=O(\e).\eeq
We have then on the boundary $$   {\partial_{yy}^2 \U_\e (\t, 0)=-\partial_{\t\t}^2 \U_\e(\t, 0)+\kappa(\t)\partial_y \U_\e (\t, 0)+\U_\e(\t, 0)}$$
and
\begin{equation*}
\begin{aligned}
   {\partial_{yyy}^3 \U_\e (\t, 0)}&=   {-2\kappa(\t)\partial_{\t\t}^2 \U_\e(\t, 0)-\partial^3_{\t\t y}\U_\e(\t, 0)-\dot\kappa(\t)\partial_\t \U_\e(\t, 0)}\\
&   {+\kappa^2(\t)\partial_y\U_\e(\t, 0)+\kappa(\t)\partial_{yy}^2 \U_\e(\t, 0)+\partial_y\U_\e(\t, 0)}
\end{aligned}
\end{equation*}
Then differentiating twice with respect to $\t$ the value $\U_\e$ on the boundary and the values of $\partial_\nu \U_\e$ on the boundary we get
\begin{equation}\label{ue2}
\begin{aligned}
&   {\partial_{\t\t}^2 \U_\e(\t, 0)}=   {-\frac{\e}{\sqrt 2}\left(\frac{d^2}{d\t^2}\ln\hat\mu^2-\e\frac{d^2}{d\t^2}(\hat\mu \nu_2(\t))-\e^2 \frac{d^2}{d\t^2}(\hat \mu^2 \zeta_2(\t))\right)= -\frac{\varepsilon}{\sqrt 2}\frac{d^2}{d\t^2}\ln\hat\mu^2+O(\e^2)}\\
&   {\partial_{\t\t y}^3 \U_\e(\t, 0)}=   {-\frac{d^2}{d\t^2}\frac{1}{\hat\mu}+\frac{\e}{\sqrt 2}\left(2\ddot\kappa(\t)+\ddot{\hat\mu}\ln 4 +\e \frac{d^2}{d\t^2}(\hat\mu \zeta_1)\right)}=   {-\frac{d^2}{d\t^2}\frac{1}{\hat\mu}+O(\e)}
\end{aligned}
\end{equation}
By \eqref{ue1} and \eqref{ue2}
we deduce
\beq\label{ue3} 
\begin{aligned}
&   {\frac{\sqrt 2}{\e}\partial^2_{yy}\mathcal U_\e(\theta,0)=\frac{\sqrt 2}{\e} \[1-\frac{\kappa(\t)}{\hat\mu}+\frac{\e}{\sqrt 2}\left(\frac{d^2}{d\t^2}\ln \hat\mu^2-\ln\hat\mu^2+2\kappa^2(\t)+\kappa(\t)\hat\mu \ln 4\right)\]+O(\e)}
 \\
&   {\frac{\sqrt 2}{\e}\partial^3_{yyy}\mathcal U_\e(\theta,0)}=   {\frac{\sqrt 2}{\e}\[-\frac1{\hat\mu}-\frac{2}{\hat\mu}\kappa^2(\t)+\kappa(\t)   {+\frac{d^2}{d\t^2}\frac{1}{\hat\mu}}\]+O(1)} .
\end{aligned}
\eeq
Finally, by \eqref{glu1}, \eqref{glu2} and \eqref{ue3} the claim follows.

\end{Proof}

 \subsection{The   main term    in the whole domain}
The main term of the solution is given by
\beq\label{baru}
U_\lambda(x) = \eta_\delta(y(x))u_{  \l}(\theta(x),y(x))+\Big(1-\eta_\delta(y(x))\Big){\sqrt2\over\e}\mathcal U_\e(x), \eeq
where $u_{  \l}$ is defined in \eqref{acb}, $\mathcal U_\e$ is defined in  \eqref{pbfuoriep1} and
 $\eta_\d(x)=\eta_\d(y(x))$ is a cut-off function  such that $\eta_\d=1$ in $\mathcal D_\d,$  $\eta_\d=0$ in $\Omega\setminus\mathcal D_{2\d},$ $0\le\eta_\d\le1$
 and $|\eta'_\d|\le \frac 1\delta$ and $|\eta^{''}_\d|\le \frac 1{\delta^2}.$
  We choose (see Lemma \eqref{starnorma})
 \beq\label{delta}
\d:=\e^a\qquad a\in \left(   {\frac{13}{14}}, 1\right).
\eeq

  \section{The error estimate}\label{2}

In this section we study the error term
\beq\label{erroreS}
  S_\lambda(U_\l):=-\Delta U_\l +U_\l-\l e^{U_\l}, \qquad \mbox{in}\,\, \O.
\eeq

 \subsection{Estimate of the error close to the boundary}
  It is useful to scale the problem. After the change of variables \eqref{cambio}, in a neighborhood of the curve, we get that the error term is given by
 \beq\label{Rv}
 \mathcal R ( \tilde U_\l):= \hat\mu^2_0(\e s)\partial_{ss}^2 \tilde U_\l +\partial^2_{tt} \tilde U_\l +\tilde{\mathcal A}( \tilde U_\l) +\l \mu^2 e^{\tilde U_\l}\qquad \mbox{in}\,\,\ \mathcal C_{2\d}
 \eeq
 where $\tilde{\mathcal A}$ is the operator defined in \eqref{tildeA} and $ \tilde U_\lambda$ is defined as follows:
\beq\label{approxcambio}
\tilde U_{\l} (s,t):=\left\{
\begin{array}{lr}
\tilde u_{  \l}(s,t) \qquad \mbox{in}\,\,\ \mathcal C_\d\\\\
\tilde\eta_\d (t)\tilde u_{  \l}(s, t)+(1-\tilde\eta_\d(t)) \frac{\sqrt 2}{\e}\U_\e( \e s, \mu t) \qquad \mbox{in}\,\,\,\, \mathcal C_{2\d}\setminus \mathcal C_\d
\end{array}
\right.
\eeq
where $\tilde u_{  \l}$ is the scaled function $u_\lambda$ defined in \eqref{acb}, i.e.
\beq\label{approxcambiov1}
\tilde u_{  \l}(s, t):= \ln\frac{4}{\e^2}-\ln\l +\ln\frac{1}{4\hat\mu(\e s)^2}+w(t)+\tilde\alpha_\mu(s, t)+\tilde v_\mu(s, t)+\tilde\beta_\mu(s, t)+\tilde z_\mu(s, t)+ e_0^\e(\e s) Z_0(t).
\eeq
 Here $\tilde\eta_\d( t)= \eta_\d(\mu t)$ is the cut-off function $\eta_\d$ scaled, which is $1$ inside $\mathcal C_\d$ and $0$ outside $\mathcal C_{2\d}$.
It is only  necessary to compute the rate of the error part $\tilde{\mathcal R}(\tilde U_\l )$ defined as
\beq\label{rtilde}\begin{aligned}\tilde{\mathcal R}(\tilde U_\l )&:= \mathcal R(\tilde U_\l) -  \tilde\eta_\d [\e^2\hat\mu_0^2 \ddot e_0^\e(\e s)+\Lambda_1 e_0^\e(\e s)]Z_0(t)\end{aligned}\eeq

 \begin{lemma}\label{starnorma}
There exist $C>0$ and $\e_0>0$ such that for all $\e\in (0, \e_0)$ we get
$$\|\tilde{\mathcal R}(\tilde U_\l)\|_{**}\leq C    {\e^{\frac 52}} .$$
\end{lemma}
 \begin{Proof}

 For sake of simplicity, let  $v_{2, \l} :=\tilde\eta_\d v_{1, \l} +(1-\tilde\eta_\d)v_{3, \l} $ where $v_{1, \l} (s, t)=\tilde u_\l(s,t)$ and
 $v_{3, \l}(s, t):=\frac{\sqrt 2}{\e}\U_\e( \e s, \mu t)$.
 We are going to estimate $\|\tilde{\mathcal R}(v_{i, \l})\|_{**}$ for $i=1,2,3.$\\\\

 It is useful to point out that  the weight $1+|t|^{2+\sigma}$ present in the weighted norm $\|\cdot\|_{**}$  in $\mathcal C_{2\d} $ has the following growth
\beq\label{peso}
\sup_{\mathcal C_{2\d} } (1+|t|^{\sigma+2})=O\(\e^{-(1-a)(\sigma+2)}\).
\eeq

 {\em Claim 1:} $\|\tilde{\mathcal R} (v_{1, \l})\|_{**}\leq C    {\e^{\frac 5 2}}$.\\\\
 For sake of simplicity, set
 $$\tilde h_\mu(s,t):=h_\mu(\e s,\mu   {(\e s)} t)\ \hbox{and}\ h_\mu(\t,y):=\alpha_\mu(\t,y)+v_\mu(\t,y)+\beta_\mu(\t,y)+z_\mu(\t,y).$$
 We have to take into account that $\mu=\e\hat\mu$.
 Therefore, a direct computation proves that
 \beq\label{uf1}
 \partial_t\tilde h_\mu =\e\hat\mu\partial_yh_\mu \ \hbox{and}\ \partial_s\tilde h_\mu(s,t)=\e \partial_\t h_\mu  +\e^2\dot{\hat\mu}  t\partial_yh_\mu,
 \eeq
 \beq\label{uf2}
 \partial_{tt}^2\tilde h_\mu =\e^2\hat\mu^2\partial_{yy}^2h_\mu ,\eeq
  \beq\label{uf3}\partial_{ss}^2\tilde h_\mu =\e^2 \partial_{\t\t}^2 h_\mu +   {2}\e^3\dot{\hat\mu}  t\partial_{\t y}^2h_\mu  +\e^3\ddot{\hat\mu}  t\partial_{  y}h_\mu +\e^4\dot{\hat\mu}^2  t^2\partial_{y y}^2h_\mu
 \eeq
 and
  \beq\label{uf4}\partial_{st}^2\tilde h_\mu =\e^2 \dot{\hat\mu}\partial_{y} h_\mu +\e^2 {\hat\mu}  \partial_{\t y}^2h_\mu  +\e^3\hat \mu\dot{\hat\mu}  t\partial_{  yy}^2h_\mu.
 \eeq

A straightforward computation  together with   Lemmas \ref{lemmaalpha}, \ref{lemmav},  \ref{lemmabeta} and \ref{lemmaz} lead to

\begin{equation}\label{rtv1}
\begin{aligned}
\tilde{\mathcal R}(v_{1, \l})&=   {\mathcal R \( \ln \frac{4}{\e^2}-\ln\l +\ln\frac{1}{4\hat\mu^2}+w+\tilde h_\mu +e_0^\e(\e s) Z_0(t)\)-[\e^2\hat\mu_0^2 \ddot e_0^\e(\e s)+\Lambda_1 e_0^\e(\e s)]Z_0(t)}\\
 &=   {\frac{\hat\mu^2}{(1-\mu t \kappa)^2}\partial_{ss}^2\tilde h_\mu}+    {\frac{\dot\mu^2 t^2}{(1-\mu t \kappa)^2}\partial_{tt}^2 \tilde h_\mu}   {-\frac{2\e^{-1}\mu\dot\mu t}{(1-\mu t \kappa)^2}\partial_{st}^2\tilde h_\mu}    {-\frac{\mu^2\dot\mu t^2 \dot\kappa}{(1-\mu t \kappa)^3}\partial_t w} \\
&-\frac{\mu\kappa}{1-\mu t \kappa}\partial_t \tilde z_\mu+\frac{\mu^3\e^{-1}t\dot\kappa}{(1-\mu t\kappa)^3}\partial_s\left[\ln\frac{1}{4\hat\mu^2}+\tilde h_\mu\right]-\mu^2\tilde h_\mu \\
&   {+\(-\frac{\ddot\mu\mu t}{(1-\mu t \kappa )^2}-\frac{\mu^2\dot\mu t^2 \dot\kappa}{(1-\mu t\kappa)^3}+\frac{2\dot\mu^2 t}{(1-\mu t\kappa)^2}\right) \partial_t \tilde{h}_\mu }+\tilde{\mathcal A}\( e^\e_0(\e s) Z_0(t)\)\\
&+\underbrace{e^w\left(e^{\tilde h_\mu+e_0^\e Z_0}-1-\tilde v_\mu -\mu\alpha_1-\tilde z_\mu -\mu^2\left(\alpha_2+\beta_1+\frac 12(\alpha_1+ v)^2\right)-e_0^\e Z_0 \right)}_{\mathcal S_0}.
\end{aligned}
\end{equation}
Now 
\begin{eqnarray}\label{buf1}
   {\tilde{\mathcal R}(v_{1, \l})-\tilde{\mathcal A}\(e^\e_0(\e s) Z_0(t)\)-\mathcal S_0}  &=&  O\(\partial_{ss}^2\tilde h_\mu \)+O\(\e^2   {|t|^2}\partial_{tt}^2\tilde h_\mu \)+O\(\e    {|t|}\partial_{st}^2\tilde h_\mu \)\nonumber\\
& &+O\(   {\e^3 |t|^2\partial_{t}  w} \)+O\(   {(\e^2 |t|+\e^3 |t|^2)\partial_{t} \tilde h_\mu }\)+O\(\e^2|t|\partial_s\tilde h_\mu \)\nonumber\\
& &+O\(\e \partial_{t} \tilde z_\mu \)+O\(\e^2   \tilde h_\mu \)+O\(\e ^3 |t|\)
\end{eqnarray}
By using the estimates of Lemma \ref{uffa} together with the derivatives of the function $\tilde h_\mu$ computed in \eqref{uf1}, \eqref{uf2}, \eqref{uf3} and \eqref{uf4},  we get
$${\| \tilde{\mathcal R}(v_{1, \l})-\tilde{\mathcal A}\(e^\e_0(\e s) Z_0(t)\)-\mathcal S_0\|_{**}\leq c\e^{4a-1-(1-a)(\sigma+2)}}$$
 Now
\begin{equation*}
\begin{aligned}
   {\tilde{\mathcal A}\(e^\e_0(\e s) Z_0(t)\)} &   { := b_0 (s, t) \e^2 \ddot e_0^\e Z_0+ b_1(s, t)  e_0^\e \partial^2_{tt} Z_0+b_2(s, t)\e \dot e_0^\e \partial_ t Z_0}\\
&   {+b_3(s, t) \left[e_0^\e \partial_t Z_0\right]+b_4(s, t) \left[ \e\dot e_0^\e Z_0\right]-\mu^2  e_0^\e Z_0}
\end{aligned}
\end{equation*}
where the coefficients $b_j(s, t)$ are defined in \eqref{tildeA}.
Then we deduce immediately
\begin{equation*}
\|\tilde{\mathcal A}\(e^\e_0(\e s) Z_0(t)\)\|_{**}\leq  c \e^{   {\frac 5 2}}
\end{equation*}
Finally,
\begin{equation*}
\begin{aligned}
   {\mathcal S_0} & \underbrace{   {= e^w\left[ e^{\tilde h_\mu} -1-\tilde v_\mu -\mu\alpha_1 -\tilde z_\mu -\mu^2 \(\alpha_2+\beta_1 +\frac 12 (\alpha_1+ v)^2\)\]}}_{\mathcal S_0^1}\\
&\underbrace{{   {+ e^{w+\tilde h_\mu }\left[e^{e_0^\e Z_0}-1-e^\e_0Z_0\right]}}}_{\mathcal S_0^2}+ \underbrace{{   {e^w \(e^{\tilde h_\mu}-1\) e_0^\e Z_0 }}}_{\mathcal S_0^3}
\end{aligned}
\end{equation*}
and hence $${\|\mathcal S_0^1\|_{**}\leq C \e^3}$$ and it is independent of $e_0$ while $${\|\mathcal S_0^2\|_{**}\leq C \e^{2\gamma_0}}$$ and it is quadratic in $e_0$ and finally $${\|\mathcal S_0^3\|_{**}\leq C \e^{\frac 5 2 }}$$ and it is linear in $e_0$. Since $   {a>\frac{13}{14}}$ the claim follows.

{\em  Claim 2:} $\|\tilde{\mathcal R} (v_{3, \l})\|_{**}\leq C e^{-\frac{c}{\e}}$ for some positive constant $c.$\\\\
Since $$\tilde{\mathcal R}(v_{3, \l})= {\mathcal R}(v_{3, \l})=\l \mu^2 e^{v_{3, \l}},$$ we get
$$\|\mathcal R(v_{3, \l})\|_{**}\leq c \sup_{\mathcal C_{2\d}\setminus \mathcal C_{\d}}\left|\frac{4}{\e^2}e^{-\frac{\sqrt 2}{\e}} e^{v_{3, \l}}\right|(1+|t|^{\sigma+2})\leq c e^{-\frac{c}{\e}}.$$

{\em Claim 3:} $\|\tilde{\mathcal R} (v_{2, \l})\|_{**}\leq C    {\e^{\frac 52}}$.\\\\

By making some tedious computations, one gets that
\beq\label{espRv2}
\tilde{\mathcal R}( v_{2, \l})= \tilde\eta_\d \tilde{\mathcal R}(v_{1, \l})+(1-\tilde\eta_\d)\mathcal R(v_{3, \l})+\mathcal R_2
\eeq
where
\begin{equation*}
\begin{aligned}
   {\mathcal R_2} &:=   { \underbrace{\hat\mu_0^2\partial_{ss}^2\tilde\eta_\d(v_{1, \l}-v_{3, \l})+2\hat\mu_0^2\partial_s \tilde\eta_\d\partial_s ( v_{1, \l}-  v_{3, \l})+\partial_{tt}^2\tilde\eta_\d ( v_{1, \l}-v_{3, \l})+2\partial_t\tilde\eta_\d \partial_t( v_{1, \l}-v_{3, \l})}_{\mathcal B_1}}\\
&+\underbrace{\tilde{\mathcal A}(\tilde\eta_\d v_{1, \l}+(1-\tilde\eta_\d)v_{3, \l})-\tilde\eta_\d\tilde{\mathcal A}(v_{1, \l})-(1-\tilde\eta_\d)\tilde{\mathcal A}(v_{3, \l})}_{\mathcal B_2}\\\\
&+\underbrace{\l \mu^2 \(e^{\tilde\eta_\d v_{1, \l}+(1-\tilde\eta_\d)v_{3, \l}}-\tilde\eta_\d e^{v_{1, \l}}-(1-\tilde\eta_\d)e^{v_{3, \l}}\)}_{\mathcal B_3}
\end{aligned}
\end{equation*}

By using the expansion \eqref{espRv2} and the result of claim 1 and claim 2 we get that
$$
\|\tilde{\mathcal R}(v_{2,\l})\|_{**}\leq c   {\e^{\frac 52}}+c e^{-\frac c\e}+\|\mathcal R_2\|_{**}
$$
So we are left to estimate $\|\mathcal R_2\|_{**}.$ Let us prove that 
\beq\label{mathcalr2}
\|\mathcal R_2\|_{**}\le C \e ^{3a -(1-a)(\sigma+2)}.
\eeq\label{mathcalr2}
Now we take into account that
$$
   {\partial_s\tilde\eta_\d=O\(\e\),\ \partial_t \tilde\eta_\d=O\( \frac\e\delta\),\ \partial^2_{ss}\tilde\eta_\d=O\(\e^2\),\ \partial^2_{st}\tilde\eta_\d=O\( \frac{\e^2}\delta\),\
\partial^2_{tt}\tilde\eta_\d=O\( \frac{\e^2}{\delta^2}\).}
$$
By \eqref{u1u3}, \eqref{u1primomenou3primo} and \eqref{uteta}, we immediately
deduce that
 for any $(s, t)\in \mathcal C_{2\d}\setminus \mathcal C_\d$ (remember that $\mu=\e\hat\mu$)
$$
 v_{1, \l}(s, t)-v_{3, \l}(s, t)=\hat u_\lambda(\e s,\mu t)-{\sqrt2\over\e}\mathcal U_\e(\e s,\mu t)=  O( \e^3 |t|^2)+O( \e^3|t|^4),
$$
$$
\partial_t \( v_{1, \l}(s, t)-  v_{3, \l}(s, t)\)=\mu \partial_y \(\hat u_\lambda(\e s,\mu t)-{\sqrt2\over\e}\mathcal U_\e(\e s,\mu t)\)=O( \e^3 |t|)+ O(   \e^3 |t|^3)
$$
and by using \eqref{uteta}
\begin{align*}
\partial_s \(v_{1, \l}(s, t)- v_{3, \l}(s, t)\)&=\e \partial_\t \(u_\lambda(\e s,\mu t)-{\sqrt2\over\e}\mathcal U_\e(\e s,\mu t)\)+ \e^2\dot{\hat \mu}t\partial_y \(u_\lambda(\e s,\mu t)-{\sqrt2\over\e}\mathcal U_\e(\e s,\mu t)\)\nonumber\\ &=O(  \e^2| t|) + O(\e^3 |t|^5)+ O(\e^4 |t|^2).
\end{align*}
Then
\begin{equation*}
\begin{aligned}
   {\mathcal B_1} &:=    {O\( \e^2 | v_{1, \l}-v_{3, \l}|\) + O\(\e |\partial_s(v_{1, \l}-v_{3, \l})|\)+O\( {\e^2 \over \d^2} | v_{1, \l}-v_{3, \l}|\)+ O\( {\e\over \d} |\partial_t ( v_{1, \l}-v_{3, \l})|\)}\\
&=    {O( \e^5|t|^4)+O(\e^3 |t|)+O\( {\e^5 \over \d^2}|t|^4\)+ O\( {  \e^4\over \d} |t|^3\)}
\end{aligned}
\end{equation*}
from which it follows that $$\|\mathcal B_1\|_{**}\leq c \e^{2a+1-(1-a)(\sigma+2)}.$$

Now
\begin{equation*}
\begin{aligned}
   {\mathcal B_2}&:=    {b_0\partial_{ss}^2 \tilde\eta_\d ( v_{1, \l}-v_{3, \l}) +
  2 b_0 \partial_{s} \tilde\eta_\d \partial_s( v_{1, \l}-v_{3, \l})+ b_1 \partial_{tt}^2 \tilde\eta_\d (	 v_{1, \l}-v_{3, \l})+2b_1 \partial_{t} \tilde\eta_\d \partial_t( v_{1, \l}-v_{3, \l})}\\
&   {+b_2 \partial_{st}^2 \tilde\eta_\d (\hat v_{1, \l}-v_{3, \l})
+b_2 \partial_{s} \tilde\eta_\d \partial_t( v_{1, \l}-v_{3, \l})+b_2 \partial_{t} \tilde\eta_\d \partial_s( v_{1, \l}-v_{3, \l})+b_2 \partial_{s} \tilde\eta_\d\partial_t ( v_{1, \l}-v_{3, \l})}\\
&   {+b_3 \partial_{t} \tilde\eta_\d ( v_{1, \l}-v_{3, \l})+b_4\partial_s \tilde\eta_\d \(v_{1, \l}-v_{3, \l}\)}
\end{aligned}
\end{equation*}
and straightforward computations  show that $$   {\|\mathcal B_2\|_{**}\le C \e^{3a-(1-a)(\sigma+2)}.}$$ Finally,
\begin{equation*}
\begin{aligned}
   {\mathcal B_3}&:=   {\l \mu^2 e^{v_{1, \l}}\(e^{(1-\tilde \eta_\d)(v_{3, \l}-v_{1, \l})}-1\)+\l\mu^2 e^{v_{1, \l}}(1-\tilde\eta_\d)\(1- e^{v_{3, \l}-v_{1, \l}}\)}
\end{aligned}
\end{equation*}
and so $$   {\mathcal B_3}=    {O\( e^w| v_{1, \l}-v_{3, \l}|\)}$$ and hence $$\|\mathcal B_3\|_{**}\le C \e^3.$$ Putting together all these estimates \eqref{mathcalr2} follows by using the fact that $a<1$. The result of the claim follows since $   {a>????}.$

That concludes the proof.

 \end{Proof}

 \begin{lemma}\label{uffa}
 Let $\alpha_\mu$, $v_\mu,$ $\beta_\mu$ and $z_\mu$ as in Lemmas \ref{lemmaalpha}, \ref{lemmav}, \ref{lemmabeta}
 and \ref{lemmaz}, respectively.
 It holds true  that uniformly  with respect to $y\in \mathcal D_{2\delta}$
 \begin{equation}\label{auffa}
\left\{\begin{aligned} &    {\alpha_\mu(\t, y)= O\(|y|^2\over \e\)+ O\(|y|^4 \over \e^2\)}\\& \partial_\t \alpha_\mu(\t,y)=O\({|y|^2\over\e}\)+O\({|y|^3\over\e^2}\)    {+ O\( |y|^5 \over \e^3\)},  \\
 & \partial_{\t\t}^2 \alpha_\mu(\t,y)=O\({|y|^2\over\e}\)+O\({|y|^3\over\e^2}\)+O\({|y|^4\over\e^3}\)    {+O\( |y|^6 \over \e^4\)}, \\
&    {\partial_y \alpha_\mu(\t, y) = O\( |y| \over \e\)},\\
&   {\partial^2_{yy}\alpha_\mu (\t, y) = O\( 1 \over \e\)+ O\(|y|^2 \over \e^2\)},\\
   & \partial_{\t y}^2 \alpha_\mu(\t,y)=O\({|y|}\),\end{aligned}\right.\end{equation}
 \begin{equation}\label{v2uffa}
\left\{\begin{aligned}
&   {v_\mu(\t, y) = O\(|y|\)},\\
 &\partial_\t v_\mu(\t,y)=O\(|y|\),  \\
 &\partial_{\t\t}^2 v_\mu(\t,y)=   {O(|y|)+O\({|y|^4\over\e^3}\)},\\
&   {\partial_y v_\mu = O\(1\),}\\
&   {\partial_{yy}^2 v_\mu = O\( |y|^2\over \e^3\)},\\
 &\partial_{\t y}^2 v_\mu(\t,y)=   {O\({|y|^3\over\e^3}\)},\end{aligned}\right.\end{equation}
 \begin{equation}\label{buffa}
\left\{\begin{aligned}&    {\beta_\mu(\t, y)= O(|y|^2),}\\
 & \partial_\t \beta_\mu(\t,y)=O\({|y|^2 }\)+   {O\({|y|^5\over\e ^3}\)} ,  \\
 & \partial_{\t\t}^2 \beta_\mu(\t,y)=O\({|y|^2 }\)+   {O\({|y|^5\over\e^3 }\)}+   {O\({|y|^6\over\e^4}\)}, \\
&   {\partial_y \beta_\mu(\t, y)= O\(|y|\),}\\
&   {\partial_{yy}^2\beta_\mu(\t, y) = O\(|y|\)}\\
   & \partial_{\t y}^2 \beta_\mu(\t,y)=O\({|y|}\),\end{aligned}\right.\end{equation}
and
 \begin{equation}\label{zuffa}
\left\{\begin{aligned}
&   {z_\mu(\t, y)= O(\e|y|)},\\ 
&\partial_\t z_\mu(\t,y)=O\(\e|y|\),  \\
 &\partial_{\t\t}^2 z_\mu(\t,y)=O(\e|y|)+   {O\({|y|^4\over\e^2}\)},\\
&   {\partial_y z_\mu(\t, y)= O\(\e\)},\\
&   {\partial^2_{yy} z_\mu(\t, y)= O\(|y|^2\over \e^2\)},\\
 &\partial_{\t y}^2 z_\mu(\t,y)=   {O\({|y|^3\over \e^2 }\)+ O\(\e\)},\end{aligned}\right.\end{equation}
\end{lemma}

\begin{Proof}

Let us remind that $\mu=\e\hat\mu(\t).$
Then
\beq\label{y1}w_\mu(y)=\ln{1\over\e^2}+\ln{1\over\hat\mu^2}+w\({y\over\e\hat\mu}\).\eeq
It is easy to check that (taking into account \eqref{mudef})
\beq\label{y11}w_\mu(y)-\ln\lambda=  {\sqrt 2\over\e }+\ln{1\over\hat\mu^2}-{\sqrt2\over\e\hat  \mu}y+O\(1\).\eeq
Moreover, some straightforward computations show  that
\beq\label{y2}\partial_\theta w_\mu(y) =-2{\dot{\hat\mu}\over\hat\mu}-{\dot{\hat\mu}\over\hat\mu^2}w'\({y\over\e\hat\mu}\){y\over\e}=O\({|y|\over\e}\),\eeq
\beq\label{y3}\partial_y w_\mu(y)= w'\({y\over\e\hat\mu}\){1\over\e\hat\mu}=O\({1\over\e}\) \eeq
and
\beq\label{y31}\partial^2_{\t y} w_\mu(y)= -{\dot{\hat\mu}\over\hat\mu^2}w'\({y\over\e\hat\mu}\){1\over\e}
-{\dot{\hat\mu}\over\hat\mu^3}w^{''}\({y\over\e\hat\mu}\){y\over\e^2}=O\({1\over\e}\) +O\({|y|\over\e^2}\),\eeq
\beq\label{y32}\partial^2_{\t\t} w_\mu(y)= -2\frac{d}{d\t}\({\dot{\hat\mu}\over\hat\mu}\)-\frac{d}{d\t}\({\dot{\hat\mu}\over\hat\mu^2}\)
w'\({y\over\e\hat\mu}\){y\over\e}+\({\dot{\hat\mu}\over\hat\mu^2}\)^2w^{''}\({y\over\e\hat\mu}\){y^2\over\e^2}=O\({|y|^2\over\e^2}\),\eeq
   {and analogously}
\begin{equation}\label{theta3}
   {\partial^3_{\t\t\t} w_\mu(y)= O\left(\frac{|y|^3}{\e^3}\right)}
\end{equation}
   {and}
\begin{equation}\label{theta4}
   {\partial^4_{\t\t\t\t} w_\mu(y)= O\left(\frac{|y|^4}{\e^4}\right)}.
\end{equation}
Moreover
 \begin{eqnarray}\label{y33}\partial^3_{\t\t y} w_\mu(y)&=&-{d\over d\t}\({\dot{\hat\mu}\over\hat\mu^2}\)
w^{''} \({y\over\e\hat\mu}\){y\over\e^2\hat\mu}+ \({\dot{\hat\mu}\over\hat\mu^2}\)^2w^{'''}\({y\over\e\hat\mu}\){y^2\over\e^3\hat\mu}\\\\
&=&O\({1\over\e}\)
+O\({|y|\over\e^2}\)+O\({|y|^2\over\e^3}\)\end{eqnarray}

Let us estimate    { $   \alpha_\mu $ and its derivatives.}\\
By    {using \eqref{y11}, \eqref{y2}, \eqref{y3}, \eqref{y31}, \eqref{y32}, \eqref{theta3}, \eqref{theta4} and \eqref{y33} and using the expression of $\alpha_\mu$ given in Lemma \ref{lemmaalpha}}
we immediately deduce the first    {three} estimates in \eqref{auffa}.
The last    {three} estimates  in \eqref{auffa} follows by the mean value theorem taking into account the initial value data in \eqref{alpha}    {and by using the equation satisfied by $\alpha_\mu$}.\\\\

Let us estimate    {$  v_\mu $ and its derivatives}.
Since
$v_\mu(\t,y)=\e\hat \mu v\(\t,{y\over\e\hat\mu}\)$ we get    {immediately} $$   {v_\mu= O\(|y|\);\qquad \partial_y v_\mu = O\(1\); \qquad \partial^2_{yy} v_\mu=O\( |y|^2\over \e^3\)}.$$    {Moreover}
$$\partial_\t v_\mu(\t,y)=\e\dot {\hat \mu} v\(\t,{y\over\e\hat\mu}\)+ \e\hat\mu\partial_\t v\(\t,{y\over\e\hat\mu}\)-{\dot{\hat\mu}\over\hat\mu}\partial_y v\(\t,{y\over\e\hat\mu}\)y=O(|y|),$$
$$\partial_{\t y}^2 v_\mu(\t,y)={\dot{\hat\mu}\over \hat\mu } \partial_y v\(\t,{y\over\e\hat\mu}\)+   \partial_{\t y}^2 v\(\t,{y\over\e\hat\mu}\)-{\dot{\hat\mu}\over\e\hat\mu^2}\partial_{yy}^2 v\(\t,{y\over\e\hat\mu}\)y-{\dot{\hat\mu}\over\hat\mu}\partial_y v\(\t,{y\over\e\hat\mu}\)=   {O\({|y|^3\over\e^3}\)},$$
and
\begin{eqnarray*}\partial_{\t\t}^2 v_\mu(\t,y)&=&\e\ddot {\hat \mu} v\(\t,{y\over\e\hat\mu}\)+ 2\e\dot{\hat\mu}\partial_\t v\(\t,{y\over\e\hat\mu}\)-{\dot{\hat\mu}^2\over\hat\mu^2}\partial_y v\(\t,{y\over\e\hat\mu}\)y\\
& &+\e {\hat \mu}   \partial_{\t\t}^2 v\(\t,{y\over\e\hat\mu}\)-2{\dot {\hat \mu}  \over\hat\mu}\partial_{\t y}^2 v\(\t,{y\over\e\hat\mu}\)y\\
& &-{d\over d\t}{\dot{\hat\mu}\over\hat\mu }\partial_{y }  v\(\t,{y\over\e\hat\mu}\)y+{\dot{\hat\mu}^2\over\e\hat\mu^3}\partial_{yy}^2 v\(\t,{y\over\e\hat\mu}\) y^2\\
& =&   {O(|y|)+O\({|y|^4\over\e^3}\)}.\end{eqnarray*}
We have used the following facts. Since $v$ solves equation \eqref{vuffa},  the functions $\partial_\t v$ and $\partial_{\t\t}^2 v$ solve  the equations
$$-\partial^2_{yy}\partial_\t v -e^w\partial_\t v=e^w\partial_\t \alpha_1(\t,y)\ \hbox{in}\ \mathbb R $$
and
$$-\partial^2_{yy}\partial_{\t\t}^2 v -e^w\partial_{\t\t}^2 v=e^w\partial_{\t\t}^2 \alpha_1(\t,y)\ \hbox{in}\ \mathbb R .$$
Therefore we apply Lemma \ref{massimo} and we deduce that $v $, $\partial_\t v$ and $\partial_{\t\t}^2 v$
 have a linear growth, namely they satisfy  for any $y\in \mathbb R$ and $\t\in[0,\ell],$ the inequalities
$$ |v (\t,y)|,|\partial_\t v (\t,y)|,|\partial_{\t\t}^2 v (\t,y)|\le c_1|y|+c_2$$
and
$$ |\partial_y v (\t,y)|,
|\partial_{\t y}^2 v(\t,y)|, |\partial_{\t\t y}^3 v   (\t,y)| \le c_3 $$
for some positive constants $c_1,c_2$ and $c_3$.
 We also remark that by equation \eqref{vuffa} we deduce that
$$|\partial_{yy}^2 v\(\t,y\)| \le a_1|y|^2+a_2|y|+a_3 \ \hbox{for any $y\in \mathbb R$ and $\t\in[0,\ell],$}$$
for some positive constants $a_1,a_2$ and $a_3$.\\\\

Arguing in a similar way, we prove estimates involving the functions $\beta_\mu$ and $z_\mu$.
\end{Proof}

 \begin{lemma}
Let $\mathcal U_\varepsilon $ be given in Lemma \ref{ue}.
Then if $\varepsilon$ is small enough
\beq\label{uteta}
 \partial_\theta\[\hat u_{ \l}(\t, y)-{\sqrt2\over\varepsilon}\U_\e(\t, y)\]= O\( |y|\)+ O\(\frac{|y|^2}{\e}\) + O\(\frac{|y|^3}{\e^2}\)+ O\(\frac{|y|^5}{\e^3}\)  \ \hbox{uniformly in}\ \mathcal D_{2\delta}
.
 \eeq
   \end{lemma}

\begin{Proof}
First of all, by mean value theorem we get for some $\bar y\in [0,y]$
$$\partial_\theta \U_\e(\t, y) =\partial_\theta \U_\e(\t, 0)+y\partial_y\( \partial_\theta \U_\e\)(\t, 0)+\frac{\sqrt 2}{\e}\partial^2_{yy}\partial_\t \U_\e (\t, \bar y) y^2 =-2\frac{\dot{\hat\mu}}{\hat\mu^2}+\frac{\dot{\hat\mu}}{\hat\mu^2} \frac{\sqrt 2}{\e} y +O\(|y|\)+ O\(\frac{|y|^2}{\e}\).$$
Here we use  the boundary condition in \eqref{pbfuoriep1} and the fact that $\partial^2_yy\( \partial_\theta \U_\e\)$ is uniformly bounded because of \eqref{ueo}.
\\
Now let us compute
\beq\label{y0}
\begin{aligned}
\partial_\theta {u}_{  \l}(\t, y) &=\partial_\theta w_\mu(y) +\partial_\theta \alpha_\mu(\t, y)+\partial_\theta v_\mu(\t, y)+\partial_\theta \beta_\mu(\t, y)+\partial_\theta z_\mu(\t, y)+ e_\e^0 Z_0\\
&= \partial_\t w_\mu (y) +O\( |y|\)+ O\(\frac{|y|^2}{\e}\) + O\(\frac{|y|^3}{\e^2}\)+ O\(\frac{|y|^5}{\e^3}\) \\
&=-2\frac{\dot{\hat\mu}}{\hat\mu^2}+\frac{\dot{\hat\mu}}{\hat\mu^2} \frac{\sqrt 2}{\e} y +O\( |y|\)+ O\(\frac{|y|^2}{\e}\) + O\(\frac{|y|^3}{\e^2}\)+ O\(\frac{|y|^5}{\e^3}\) .
\end{aligned}
\eeq
We take into account estimate \eqref{y2} together with the first estimates in \eqref{auffa},  \eqref{v2uffa},  \eqref{buffa}
and  \eqref{zuffa}. 
 Then the claim follows.

\end{Proof}

  \subsection{Estimate of the error in the inner part  }

  \begin{lemma}\label{errore1}
 There exist $c>0$ and $\e_0>0$ such that for any $\e\in(0, \e_0)$ we have
 \beq
 \|  S_\l(U_\l)\|_{\infty,\Omega\setminus D_{2\delta}}\leq e^{-{c\over\e}}.
 \eeq
\end{lemma}
\begin{Proof}

Since $\mathcal U_\e$ solves \eqref{pbfuoriep1} we have
$$S_\l(U_\l) =S_\l\({\sqrt2\over\e}\mathcal U_\e\)= -\l e^{\frac{\sqrt{2}}{\e}\U_\e}=-\frac{4}{\e^2}e^{\frac{\sqrt{2}}{\e}(\U_\e-1)}
=-\frac{4}{\e^2}e^{\frac{\sqrt{2}}{\e}(\U_\e-\U_0)}e^{\frac{\sqrt{2}}{\e}(\U_0-1)}.
$$
Now, by the fact that  $\partial_\nu \U_0 <0$ we   deduce that $\U_0(x)-1\le c<0$ if $x\in\O\setminus \mathcal D_{2\delta}$ for some constant $c.$
Moreover, by \eqref{ueo} we also deduce that
$|\U_\e(x)-\U_0(x)|\le c\varepsilon$ for any $x\in\O\setminus \mathcal D_{2\delta}$ for some constant $c.$
 Therefore, the claim follows. \end{Proof}

 \subsection{The projection of the error   along  ${Z_0}$ }\label{proiezioni} 
We are going to  compute the component of the scaled error $\mathcal R (\tilde U_\l)$ given in \eqref{Rv} along $Z_0$.
\begin{lemma}\label{lemmaproiezioni}
There exists $\e_0>0$ such that for any $\e\in (0, \e_0)$ the following expansion hold:
\begin{equation}\label{proiezionez1}
\begin{aligned}
\int_{-\frac{2\d}{\mu}}^0 \mathcal R ( \tilde U_\l) Z_0\, dt &= \e^{   {\frac 32}}\left[\e^2 \(a_0(\e s)\ddot e_0(\e s)+ a_1^\e(\e s) \dot e_0\)+ a_2(\e s) e_0\]+ \e^3 M_0(\e s) \\
&+\e^{   {3}}H_0(e_0, \dot e_0, \ddot e_0)  \ \hbox{for any $s\in [0, \frac{\ell}{\e}]$, }
\end{aligned}
\end{equation}
where
\beq\label{a0}
a_0(\e s)= \hat\mu^2_0 + \e  a^\e_0(\e s)
\eeq

and 
\beq\label{a2}
a_2(\e s)= \Lambda_1 + \e  a_2^\e (\e s)
\eeq
with $a^\e_i$ $i=0, 1,2$ explicit smooth functions, uniformly bounded in $\e$. Moreover in \eqref{proiezionez1}
\begin{itemize}
\item $M_0$ is a sum of  explicit smooth functions of the form, uniformly bounded in $\e$;
\item $H_0$ denotes a sum of functions of the form
$$h_0(\e s)\[ h_1(e_0)+o(1) h_2(e_0, \dot e_0, \ddot e_0)\]$$

\begin{itemize}
\item $h_0$ is a smooth function uniformly bounded in $\e$;
\item $h_1$ and $h_2$ is a smooth function of its arguments, uniformly bounded in $\e$ when $e_0$ satisfys \eqref{normae};

\item $o(1)\rightarrow 0$ as $\e\rightarrow 0$ uniformly when $e_0$ satisfys \eqref{normae}.
\end{itemize}
\end{itemize}
\end{lemma}
\begin{Proof}
For sake of simplicity, let  $v_{2, \l} :=\tilde\eta_\d v_{1, \l} +(1-\tilde\eta_\d)v_{3, \l} $ where $v_{1, \l} (s, t)=\tilde u_\l(s,t)$ and
 $v_{3, \l}(s, t):=\frac{\sqrt 2}{\e}\U_\e( \e s, \mu t)$.
\\
First of all we get that by using \eqref{rtilde} and \eqref{espRv2}
\begin{eqnarray*}
\int_{-\frac{2\d}{\mu}}^0 \mathcal R (\tilde U_\l)Z_0\, dt &=& \underbrace{\int_{-\frac{2\d}{\mu}}^{-\frac{\d}{\mu}} \mathcal R_2 Z_0\, dt}_{I_1^0}+\underbrace{\int_{-\frac{2\d}{\mu}}^{-\frac{\d}{\mu}} \mathcal R(v_{3, \l}) Z_0\, dt}_{I_2^0}+\underbrace{\int_{-\frac{2\d}{\mu}}^{0} \tilde\eta_\d\tilde{\mathcal R }(v_{1, \l}) Z_0(t)\, dt}_{I_3^0}\\
&&+\underbrace{\[\e^2\hat\mu_0^2 \ddot e^\e_0(\e s)+\Lambda_1  e^\e_0(\e s)\]\int_{-\frac{2\d}{\mu}}^0\tilde\eta_\d Z_0^2(t) \, dt}_{I_4^0}.
\end{eqnarray*}
We remark that

$$\int_{-\frac{\d}{\mu}}^0 Z_0^2(t)\, dt= \frac 1 2+ O\(e^{-\sqrt{\Lambda_1 }\frac{\d}{\mu}}\)$$ and hence
$$I_4^0 := \frac 12\e^{\frac 32} \[\e^2\hat\mu_0^2 \ddot e_0(\e s)+\Lambda_1  e_0(\e s)\]+ O\(e^{-\sqrt{\Lambda_1 }\frac{\d}{\mu}}\).$$ Moreover by using Claim 2 of Lemma \ref{starnorma} we get that
$$I_2^0= O\( e^{-c\frac{1}{\e}}\)$$ for some positive $c$ and similarly, by using Claim 3 of Lemma \ref{starnorma} and also the exponential decay of $Z_0$ it follows that 
$$I_1^0 =O\( e^{-c\frac{1}{\e}}\)$$ for some positive $c$. It remains to evaluate only $I_3^0$. By using  \eqref{rtv1}

\begin{equation*}
\begin{aligned}
I_3^0&= \int_{-\frac{\d}{\mu}}^{0} \tilde{\mathcal R }(v_{1, \l}) Z_0(t)\, dt+O\( e^{-{\frac c \e}}\)\\
&= \int_{-\frac{\d}{\mu}}^{0}\frac{\hat\mu^2}{(1-\mu t \kappa)^2}\partial_{ss}^2\tilde h_\mu Z_0(t)\, dt+ \int_{-\frac{\d}{\mu}}^{0}\frac{\dot\mu^2 t^2}{(1-\mu t \kappa)^2}\partial_{tt}^2 \tilde h_\mu Z_0(t)\, dt-\int_{-\frac{\d}{\mu}}^{0}\frac{2\e^{-1}\mu\dot\mu t}{(1-\mu t \kappa)^2}\partial_{st}^2\tilde h_\mu Z_0(t)\, dt\\
& - \int_{-\frac{\d}{\mu}}^{0}\frac{\mu^2\dot\mu t^2 \dot\kappa}{(1-\mu t \kappa)^3}\partial_t w Z_0(t)\, dt- \int_{-\frac{\d}{\mu}}^{0}\frac{\mu\kappa}{1-\mu t \kappa}\partial_t \tilde z_\mu Z_0(t)\, dt+\int_{-\frac{\d}{\mu}}^{0}\frac{\mu^3\e^{-1}t\dot\kappa}{(1-\mu t\kappa)^3}\partial_s\left[\ln\frac{1}{4\hat\mu^2}+\tilde h_\mu\right]Z_0(t)\, dt\\
&-\int_{-\frac{\d}{\mu}}^{0}\mu^2\tilde h_\mu Z_0(t)\, dt+\int_{-\frac{\d}{\mu}}^{0}\(-\frac{\ddot\mu\mu t}{(1-\mu t \kappa )^2}-\frac{\mu^2\dot\mu t^2 \dot\kappa}{(1-\mu t\kappa)^3}+\frac{2\dot\mu^2 t}{(1-\mu t\kappa)^2}\right) \partial_t \tilde{h}_\mu Z_0(t)\, dt\\
&+\int_{-\frac{\d}{\mu}}^{0}\mathcal S_0^1 Z_0(t)\, dt+\int_{-\frac{\d}{\mu}}^{0}\tilde{\mathcal A}\( e^\e_0(\e s) Z_0(t)\)Z_0(t)\, dt+\int_{-\frac{\d}{\mu}}^{0}\mathcal S_0^2 Z_0(t)\, dt+\int_{-\frac{\d}{\mu}}^{0}\mathcal S_0^3 Z_0(t)\, dt+O\( e^{-{\frac c \e}}\)\\
&= \e^{3}M_0(\e s)(1+o(1))+\int_{-\frac{\d}{\mu}}^{0}\tilde{\mathcal A}\( e^\e_0(\e s) Z_0(t)\)Z_0(t)\, dt+\int_{-\frac{\d}{\mu}}^{0}\mathcal S_0^2 Z_0(t)\, dt\\
&+\int_{-\frac{\d}{\mu}}^{0}\mathcal S_0^3 Z_0(t)\, dt+O\( e^{-{\frac c \e}}\)
\end{aligned}
\end{equation*}
where $M_0 (\e s)$ is a sum of smooth functions uniformly bounded in $\e$ that does not depend on $e_0$.
Now
\begin{equation*}
\begin{aligned}
\int_{-\frac{\d}{\mu}}^{0}\tilde{\mathcal A}\( e^\e_0(\e s) Z_0(t)\)Z_0(t)\, dt&=\e^{\frac 52}\[\e^2 \ddot e_0 \underbrace{\left(\frac 12 \frac{\partial \hat\mu_\e}{\partial \e}_{|_{\e=0}}+2\kappa(\e s)\hat\mu^2(\e s) \int_{-\infty}^{0} t Z_0^2(t)\, dt\right)}_{a_0^\e (\e s)}\]\\
&+\e^{\frac 52}\[\e \dot e_0 \underbrace{\left(-2\hat\mu\dot{\hat\mu}\int_{-\infty}^0 t \partial_t Z_0 Z_0\, dt \right)}_{a_1^\e (\e s)}\]+\e^{\frac 52}\left[ e_0 \underbrace{\left(-\hat\mu \kappa \int_{-\infty}^0 \partial_t Z_0 Z_0\, dt\right)}_{\hat a_2^\e(\e s)}\right]\\
&+ \e^{\frac 7 2} h(e_0, \dot e_0, \ddot e_0)(1+o(1))
\end{aligned}
\end{equation*}
where $h(e_0, \dot e_0, \ddot e_0)$ is a sum of functions depending linearly on $e_0, \dot e_0, \ddot e_0$.
Now $$\int_{-\frac{\d}{\mu}}^{0}\mathcal S_0^2 Z_0\, dt = \e^{3} F(e_0)(1+o(1));\qquad \int_{-\frac{\d}{\mu}}^{0}\mathcal S_0^3 Z_0\, dt = \e^{\frac 32}\[ e_0 \underbrace{\int_{-\infty}^0 e^w Z_0^2 (\alpha_1+v)\, dt}_{\tilde a_2^\e (\e s)} \](1+o(1))$$
where $F$ is quadratic in $e_0$.\\ Putting together all these estimates we get
\begin{equation*}
\begin{aligned}
\int_{-\frac{2\d}{\mu}}^0 \mathcal R (\tilde U_\l)Z_0\, dt &=\e^{\frac 32}\left[\e^2 \underbrace{\left(\frac 12 \hat\mu_0^2 +\e a_0^\e \right)}_{a_0(\e s)}\ddot e_0 + a_1^\e (\e s) \dot e_0+\underbrace{\left(\frac 12\Lambda_1+\e a_2^\e \right)}_{a_2(\e s)}e_0\right]+\e^3 M_0(\e s)\\
&+\e^{3} F(e_0)(1+o(1))+\e^{\frac 72}h(e_0, \dot e_0, \ddot e_0)(1+o(1))
\end{aligned}
\end{equation*}
and the result follows.
\end{Proof}

\section{The remainder term}\label{3}

 We split the remainder term $ \Phi_\m $ in \eqref{pro1}  as
\beq\label{solgluing} \Phi_\m =\eta_{2\d}     \phi_\m+    \psi_\m,\eeq
where $    \phi_\m$ solves a linear problem defined in a neighborhood of the boundary
and $    \psi_\m$ solves a linear problem defined in the whole domain.
More precisely, we are led to consider the couple of linear problems
\begin{equation}\label{gluing1}
\left\{
\begin{aligned}
  &\Delta    \psi -    \psi+(1-\eta_{2\d})\m e^{U_\m}    \psi = -(1-\eta_{2\d}) S_\m(U_\m)-(1-\eta_{2\d})N(\eta_{2\delta}    \phi +    \psi )\\
 &\hskip5truecm -2\nabla\eta_{2\delta}\nabla    \phi -    \phi \Delta\eta_{2\delta}\qquad   \mbox{in}\,\,
\Omega\\
 &\partial_\nu     \psi  =0\  \mbox{on}\,\, \partial\Omega\\
\end{aligned}
\right.
\end{equation}
and
\begin{equation}\label{gluing2}
\left\{
\begin{aligned}
&L(    \phi ) =-S_\m(U_\m)- N(\eta_{2\delta}    \phi +    \psi )-\m e^{U_\m}    \psi  \  \mbox{in}\,\,
\mathcal{D}_{2\d}\\
&\partial_\nu     \phi   =0 \qquad \mbox{on}\,\,\, \partial\mathcal D_{2\delta}.  \\
 \end{aligned}
\right.
\end{equation}

\subsection{The remainder term in the whole domain}\label{gluing}

Given a function $    \phi$ defined in a neighborhood of the boundary, let us find a function $    \psi $ which solves problem \eqref{gluing1}.\\

First of all, it is useful to point out that
for any $g\in L^\infty(\O)$ there exists a unique $    \psi$ solution to the linear problem
\begin{equation}\label{gluing3}
\left\{
\begin{array}{lr}
\Delta    \psi -    \psi+\(1-\eta_{2\d}\) \m e^{U_\m}    \psi =g\qquad \mbox{in}\,\, \Omega\\
\partial_\nu     \psi =0\qquad\qquad\quad\qquad\mbox{on}\,\,\, \partial\Omega
\end{array}
\right.
\end{equation}
 with
 \beq\label{stimapsiprimo}
\|    \psi \|_\infty\leq C\|g\|_\infty.
\eeq
It is enough to show that  the linear perturbation term $\(1-\eta_{2\d}\) \m e^{U_\m}    \psi $ is  small   as $\e$ goes to zero. Indeed,   arguing as in Lemma \ref{errore1}, we
have
  $$\|(1-\eta_{2\d})\lambda e^{u_\l}\|_\infty\leq e^{-\frac c\e}
 $$  for some positive constant $c.$\\\\

 Now, let us split the remainder $    \psi=    \psi^1+    \psi^2$ where $    \psi^1$ solves a linear problem and $    \psi^2$ solves a nonlinear problem. More precisely, $    \psi^1$ solves \eqref{gluing3} with
\begin{equation}\label{gluing4}
g:=-(1-\eta_{2\d})S_\m(
U_\m)-(1-\eta_{2\d})N(\eta_{2\delta}    \phi )-2\nabla\eta_{2\delta}\nabla    \phi -    \phi \Delta\eta_{2\delta}
\end{equation}
and $    \psi^2$ solves \eqref{gluing3} with
\begin{equation}\label{gluing5}
g:=-(1-\eta_{2\d})[N(\eta_{2\delta}    \phi +    \psi^1 +    \psi^2 )-N(\eta_{2\delta}    \phi )].\end{equation}

It is clear that for any function $     \phi$ there exists a unique $    \psi^1$ solution to \eqref{gluing3} with the R.H.S. as in \eqref{gluing4}. Let us prove that
\beq\label{stimapsi1}
\|\psi^1_\m\|_\infty \leq c \e^{(\sigma+2)(1-a)}\|\tilde\phi\|_*.
\eeq
By \eqref{stimapsiprimo} we need to estimate the $L^\infty$-norm of R.H.S. given in \eqref{gluing4}.
First of all,  in Lemma \ref{errore1} we have
  $$\|(1-\eta_{2\d})S _\l (U_\l)\|_\infty\leq e^{-\frac c\e}
 $$  for some positive constant $c.$
 Moreover
\begin{equation*}
\begin{aligned}
 \|(1-\eta_{2\d}) N(\eta_{2\d}    \phi)\|_\infty&\leq c\|(1-\eta_{2\d})\m e^{U_\m}\eta_{2\d}^2    \phi^2\|_\infty\leq c\|(1-\eta_{2\d})\m e^{u_{ \l}}\|_\infty\|(1-\eta_{2\d})\eta_{2\d}    \phi \|_\infty^2\\
&\leq c e^{-\frac c\e}\e^{2\sigma(1-a)}\|\tilde{    \phi}\|_*
  \end{aligned}
\end{equation*}
  since $$\|(1-\eta_{2\d})\eta_{2\d}    \phi \|_\infty \leq c \(\sup_{\mathcal C_{2\d}\setminus \mathcal C_{\d}}\frac{1}{1+|t|^\sigma}\)\|\tilde{    \phi }\|_*\leq c \e^{\sigma(1-a)}\|\tilde{    \phi}\|_*,$$
  where we agree that $\tilde{     \phi}$ is nothing but the scaled function $     \phi(\e s,\mu t).$
Finally
$$\|\nabla \eta_{2\delta}\nabla    \phi \|_{\infty}\leq \frac{\e}{\d}\(\frac{\e}{\d}\)^{\sigma+1}\|\tilde{     \phi} \|_*\leq c  \e^{(\sigma+2) (1-a)}\|\tilde{     \phi}\|_*$$ and
$$\|    \phi \Delta\eta_{2\delta}\|_{\infty}\leq \frac{\e^2}{\d^2}\(\frac{\e}{\d}\)^{\sigma }\|   \tilde \phi \|_ *\leq c\e^{(\sigma+2)(1-a)}\|\tilde{     \phi}\|_*.$$ Moreover it is possible to show thta the nonlinear operator $\psi_\l$ is Lipschitz such that $$\|\psi_\l(\phi_{\l, 1}-\psi_\l(\phi_{\l, 2})\|_\infty \le c \e^{(\sigma+2)(1-a)}\|\tilde\phi_{1}-\tilde\phi_2\|_*.$$
Once we have found the function $    \psi^1,$ we   solve equation \eqref{gluing3} with R.H.S. \eqref{gluing5}. A simple contraction mapping argument (the nonlinear term $N$ is  quadratic)    yields the existence of a function
$    \psi^2 $ such that
\beq\label{stimapsi2}
\|    \psi^2 \|_\infty \leq \|(1-\eta_{2\d})\m e^{u_\l}\|_{\infty } \|     \psi^1 \|_\infty\le e ^{-\frac c\e}\|     \psi^1 \|_\infty
\eeq for some positive constant $c.$

\subsection{The  remainder term close to the boundary: a nonlinear projected problem}
In order to solve problem \eqref{gluing2}, it is necessary to   solve a nonlinear projected problem naturally associate with it. Since it is defined in a neighborhood of the boundary, it is useful to scale it.
Then we are led to study the problem: {\em given $\mu$ which satisfies \eqref{mudef} and $e_0$  which satisfies
 \eqref{normae},  find a
function  $c_0(s)$ and a function $\tilde\phi$ so that}
\begin{equation}\label{gluing21}
\left\{
\begin{aligned}
&\mathcal L( \tilde   \phi ) =-\mathcal R (\tilde U_\l)- \mathcal N_1 (\tilde\phi )+c_0
Z_0   \  \mbox{in}\,\,
\mathcal{C}_{2\d},\\
&\partial_t   \tilde\phi   =0 ,\qquad \mbox{on}\,\,\, \partial\mathcal C_{2\delta}\cap \{t=0\}, \\
&\tilde\phi\(s+\frac\ell\e,t\)=\tilde\phi\(s ,t\)\\
&\int\limits^0_{-2\frac\delta\mu}\tilde\phi\(s ,t\)Z_0(t)dt=0,
 \end{aligned}
\right.
\end{equation}
where  $\mathcal L$ is defined in \eqref{operatoreL}, $\mathcal R (\tilde U_\l)$ is defined in \eqref{approxcambio} and the superlinear term $\mathcal N_1(\tilde\phi )$ is defined by
\begin{equation}\label{pbriscalatoN}
\mathcal N_1(\tilde\phi )=\l \mu^2 e^{\tilde U_\l}\left[e^{\tilde\eta_{2\d} \tilde\phi+\tilde\psi (\phi )}-1-(\tilde\eta_{2\d}\tilde\phi +\tilde\psi (\phi ))\right]+\l\mu^2 e^{\tilde U_\l}\tilde\psi (\phi )+\(\l\mu^2e^{\tilde U_\l}-e^w\)\tilde\phi .
\end{equation}
Here $\tilde\psi (\phi )$ is the scaled function $\[ \psi (\phi )\](\e s,\mu t) $ and $\psi (\phi )$ is the solution to the problem \eqref{gluing1}.
 In \eqref{gluing21} the terms which contains $e_1^\e$ and $e_2^\e$ in $\mathcal R(\tilde U_\l)$ are encode in the last sum (see \eqref{rtilde}).\\

 By Proposition \ref{inv} $\mathcal
L $ is invertible. Hence solving \eqref{gluing21} together with
boundary, the periodic and orthogonality conditions reduces to solve a fixed point
problem, namely
\begin{equation}\label{non4}
  \tilde   \phi =\mathcal T(- \tilde{\mathcal R}( \tilde U_\l)-
\mathcal N_1(  \tilde   \phi ))=\mathcal M(  \tilde   \phi )
\end{equation}
 where $\mathcal T$ is the operator defined in Proposition \ref{inv}. \\ We will prove the following result.
\begin{proposition}\label{nonlocale}
There exist  $c>0$ and $\m_0>0$ such that for all
$\m\in(0,\m_0)$ and for any $e_0$ satisfying \eqref{normae}, the
problem \eqref{gluing21} has a unique solution $  \tilde   \phi =\tilde   \phi (e_0)$ and
$c_0=c_0(e_0)$, which satisfy
\begin{equation}\label{non3}
\|\tilde   \phi \|_{*}\leq c \e^{\frac 52}.
\end{equation}
\end{proposition}
\begin{Proof}
Let us consider the set $$\mathcal E:=\left\{  \tilde   \phi \,\,:\,\, \|  \tilde   \phi \|_{*}\leq c\e^{\frac 52}\right\}$$ for a certain
positive constant $c$. We first show that $\mathcal M$ maps
$\mathcal E$ into itself.\\ Let $  \tilde   \phi \in \mathcal E$. Then by
using Lemma \ref{starnorma}
$$\|\mathcal M(  \tilde   \phi )\|_{*}\leq C\| \tilde{\mathcal R}( \tilde U_\l)+\mathcal N_1(  \tilde   \phi )\|_{**}\leq c\e^{\frac 52}+c\|{\mathcal N}_1(  \tilde   \phi )\|_{**}.$$
We evaluate $\|{\mathcal N}_1(  \tilde   \phi )\|_{**}$.\\
\begin{eqnarray*}
\|{\mathcal N}_1(  \tilde   \phi )\|_{**}&\leq & \|\l\mu^2 e^{\tilde U_\l}(\tilde\eta_\d   \tilde   \phi +\tilde\psi (\phi ))^2\|_{**}+\|\l\mu^2 e^{\tilde U_\l}\tilde\psi_\l(\phi)\|_{**}+\|(\l\mu^2 e^{\tilde U_\l}-e^w)\tilde   \phi \|_{**}\\
&\leq & \|e^w   \tilde   \phi ^2\|_{**}+\|e^w\tilde \psi^2_\l(\phi)\|_{**}+\|e^w   \tilde   \phi \tilde \psi(\phi)\|_{**}+\|e^w \tilde\psi (\phi )\|_{**}+\|(\l\mu^2 e^{\tilde U_\l}-e^w)  \tilde   \phi \|_{**}
\end{eqnarray*}
Now
$$\|e^w  \tilde   \phi ^2\|_{**}\leq \|  \tilde   \phi \|_*^2\sup_{\mathcal C_{2\d}}e^w \frac{1+|t|^{\sigma+2}}{(1+|t|^\sigma)^2}\leq \|  \tilde   \phi \|^2_*$$
$$\|e^w \tilde\psi^2 (\phi)\|_{**}\leq \|\tilde\psi (\phi )\|_{\infty}^2\sup_{\mathcal C_{2\d}}e^w(1+|t|^{\sigma+2})\leq \e^{2(\sigma+2)(1-a)}\|\tilde\phi\|_*^2$$
$$\|e^w   \tilde   \phi \tilde\psi (\phi )\|_{**}\leq \|\tilde\psi (\phi )\|_{\infty}\|\tilde   \phi \|_*\sup_{\mathcal C_{2\d}}e^w\frac{1+|t|^{\sigma+2}}{1+|t|^\sigma}\leq \e^{(\sigma+2)(1-a)}\|\tilde   \phi \|_*^2$$
analogously
$$\|e^w \tilde\psi (\phi )\|_{**}\leq  \e^{(\sigma+2)(1-a)}\|\tilde   \phi \|_*$$
and
finally
$$\|(\l\mu^2 e^{\tilde U_\l}-e^w)  \tilde   \phi \|_{**}\leq \|e^w (\tilde\alpha_\mu+\tilde\beta_\mu+\tilde v_\mu+\tilde z_\mu+e_0^\e Z_0)  \tilde   \phi \|_{**}\leq \e^{\min\{1, \gamma_0\}} \|  \tilde   \phi \|_*.$$
Putting together all these computations we find that
\begin{equation}\label{condstarstar}
\|{\mathcal N}_1(  \tilde   \phi )\|_{**}\leq \e^{(\sigma+2)(1-a)}\|  \tilde   \phi \|_*
\end{equation}
 and the first claim is proved.\\

 We next prove that $\mathcal M$ is a contraction,
so that the fixed point problem \eqref{non4} can be uniquely solved
in $\mathcal E$.\\ Indeed, for any $\tilde\phi_{  1}$, $\tilde\phi_{  2}\in \mathcal
M$ we get (setting $\tilde\psi_1:=\tilde\psi (\phi_{ 1})$ and $\tilde\psi_2:=\tilde\psi (\phi_{ 2})$)
\begin{eqnarray*}
\|\mathcal M(\tilde\phi_{ 1})-\mathcal M(\tilde\phi_{ 2})\|_{*}&\leq & C \|{\mathcal N}_1(\tilde\phi_{  1})-{\mathcal N}_1(\tilde\phi_{  2})\|_{**}\\
&\leq & \|e^w \cdot e^{\tilde\eta_{2\d} \tilde\phi_{ 1}+\tilde\psi_1}\left(1-e^{\tilde\eta_{2\d} (\tilde\phi_{ 2}-\tilde\phi_{  1})+\tilde\psi_2-\tilde\psi_1}+(\tilde\eta_{2\d }(\tilde\phi_{ 2}-\tilde\phi_{ 1})+\tilde\psi_2-\tilde\psi_1)\right)\|_{**}\\
&&+\|e^w \left(\tilde\eta_{2\d}(\tilde\phi_{ 2}-\tilde\phi_{ 1})+\tilde\psi_2-\tilde\psi_1\right)(e^{\tilde\eta_{2\d}\tilde\phi_{ 1}+\tilde\psi_1}-1)\|_{**}\\
&&+\|e^w(\tilde\psi_1-\tilde\psi_2)\|_{**}+\|(\l\mu^2 e^{\tilde U_\l}-e^w)(\tilde\phi_{2 }-\tilde\phi_{1 })\|_{**}\\
&\leq & \|e^w \left(\tilde\eta_{2\d }(\tilde\phi_{ 2}-\tilde\phi_{1})+(\tilde\psi_2-\tilde\psi_1)\right)^2\|_{**}\\
&&+\|e^w \left(\tilde\eta_{2\d} (\tilde\phi_ {2}-\tilde\phi_{1})+\tilde\psi_2-\tilde\psi_1\right)(\tilde\eta_{2\d}\tilde\phi_{1}+\tilde\psi_1)\|_{**}\\
&&+\|e^w(\tilde\psi_1-\tilde\psi_2)\|_{**}+\|(\l\mu^2 e^{\tilde U_\l}-e^w)(\tilde\phi_{2}-\tilde\phi_{1})\|_{**}\\
&\leq & \|e^w(\tilde\phi_{2}-\tilde\phi_{1})^2\|_{**}+\|e^w(\tilde\psi_2-\tilde\psi_1)^2\|_{**} +\|e^w(\tilde\phi_{2}-\tilde\phi_{1})(\tilde\psi_2-\tilde\psi_1)\|_{**}\\
&&+\|e^w(\tilde\phi_{2}-\tilde\phi_{1})(\tilde\eta_{2\d}\tilde\phi_{1}+\tilde\psi_1)\|_{**}+\|e^w(\tilde\psi_2-\tilde\psi_1)(\tilde\eta_{2\d}\tilde\phi_{1}+\tilde\psi_1)\|_{**}\\
&& +\|e^w (\tilde\psi_1-\tilde\psi_2)\|_{**} +\|(\l\mu^2 e^{\tilde U_\l}-e^w)(\tilde\phi_{2}-\tilde\phi_{1})\|_{**}\\
&\leq & \max_{j=1, 2}\|\tilde\phi_{j}\|_* \| \tilde\phi_{2}-\tilde\phi_{1} \|_{*}+\max_{j=1, 2}\|\tilde\psi_j\|_\infty\|\tilde\psi_2-\tilde\psi_1\|_\infty\\
&&+\max_{j=1, 2}\|\tilde\psi_j\|_\infty\|\tilde\phi_{2}-\tilde\phi_{1}\|_*+\|\tilde\phi_{1}\|_*\|\tilde\phi_2-\tilde\phi_1\|_\infty+\|\tilde\psi_1\|_\infty \|\tilde\phi_2-\tilde\phi_1\|_*\\
&&+\|\tilde\phi_1\|_*\|\tilde\psi_2-\tilde\psi_1\|_\infty +\|\tilde\psi_1\|_\infty\|\tilde\psi_2-\tilde\psi_1\|_\infty\\
&\leq & \e^{\frac 5 2} \| \tilde\phi_{2}-\tilde\phi_{1} \|_{*}.
\end{eqnarray*}

 Hence $\mathcal M$ is a contraction and the proof is complete.
\end{Proof}
 \subsection{ Proof of Theorem \ref{principale} completed}\label{fap}
It only remains to find the function $e_0$ to get   the coefficient $c_0$ in \eqref{gluing21} identically equal to zero. To do this, we multiply equation \eqref{gluing21} by $Z_0$ and we integrate in $t$. Thus the equation $$c_0(\t) =0\ \hbox{for any}\ \t\in\[0,\ell\] \quad \hbox{(here $\e s=\t$)}$$ is equivalent to
\begin{equation}\label{eqdaris}
\int_{-\frac{2\d}{\mu}}^0\[ \mathcal R( \tilde U_\l) +\mathcal L(  \tilde   \phi )+\mathcal N_1(  \tilde   \phi )\]Z_0\, dt=0\ \hbox{for any}\ s\in\[0,\frac\ell\e\].
\end{equation}
We first remark that, by using \eqref{condstarstar}, it follows that 
\begin{equation}\label{defr}\int_{-\frac{2\d}{\mu}}^0\[\mathcal L(  \tilde   \phi )+\mathcal N_1(  \tilde   \phi )\]Z_0\, dt=\e^{(\sigma+2)(1-a)+\frac 5 2} r\end{equation}
   {where $r$ is the sum of functions of the form $$h_0(\e s) \[h_1( e_0, \dot e_0)+o(1) h_2(e_0, \dot e_0, \ddot e_0)\]$$ where $h_0$ is a smooth function uniformly bounded in $\e$, $h_1$ depends smoothly on $e_0$ and on $\dot e_0$ and it is bounded in the sense that $$\|h_1\|_\infty \le c \|e_0\|_\e$$ and it is compact, as a direct application of Ascoli-Arzel\'a Theorem shows. \\ The function $h_2$ depends on $e_0, \dot e_0, \ddot e_0$ and it depends linearly on $\ddot e_0$ and it is Lipschitz with $$\|h_2(e_0^1)-h_2(e_0^2)\|_\infty \le o(1)\|e_0^1-e_0^2\|_\e .$$}
By using \eqref{proiezionez1} it follows that \eqref{eqdaris} is equivalent to the following ODE
\begin{equation}\label{sistemaode}
\begin{aligned}
\e^2 \(a_0(\e s) \ddot e_0 +a_1(\e s) \dot e_0\)+a_2(\e s) e_0&=\e^{\frac 32}M_0(\e s)+\e^{\frac 32} H_0+\e^{(\sigma+2)(1-a)+1} r
\end{aligned}
\end{equation}
where $a_i(\ e s)$, $i=0, 1, 2$ , $M_0$, $F_0$ and $H_0$ are as in Lemma \ref{lemmaproiezioni} and $r$ is as in \eqref{defr}.
%because of Proposition \ref{nonlocale} and  Lemma \ref{proiezioni}.
%Here
%\begin{itemize}
%\item the functions $M_1$ and $M_2$ are explicit functions of $\t$, smooth and uniformly bounded in $\t$ given in Lemma \ref{proiezioni};
%\item the operators $H_i=H_i(e_1, e_2)$ can be decomposed as $H_i(e_1, e_2)=A_i(e_1, e_2)+  K_i( e_1, e_2)$ where
%\begin{itemize}
%\item $ K$ is uniformly bounded if $\|e_1\|_\e+\|e_2\|_\e$ is uniformly bounded and it is compact;
%\item $A_i$ depends on $(e_1, e_2),$ $(\dot e_1,\dot e_2) $ and $(\ddot e_1, \ddot e_2),$ the dependence on $(\ddot e_1, \ddot e_2)$ being linear and it satisfies
%$$\|A_i (e_1^1, e_2^1)-A_i(e_1^2, e_2^2)\|\leq o(1) \(\| e_1^1-e_1^2\|_\e+\| e_2^1-e_2^2\|_\e\)$$ uniformly for $(e_1, e_2) $ if $\|e_1\|_\e+\|e_2\|_\e$ is uniformly bounded
%\end{itemize}
%\end{itemize}
Our goal is to find a smooth periodic function $e_0$ which solves \eqref{sistemaode}. \\ In order to do this we introduce an auxiliary problem.\\ Suppose that $p_0(\t)$ is a positive $C^2(0, \ell)$ function, $p_1(\t)$ is a $C^2(0,\ell)$ function and $\e>0$ is a parameter small enough.
\\ Given an arbitrary function $f\in C^0(0, \ell)$ let us consider the problem 
\begin{equation}\label{pbaux}
\left\{\begin{aligned}
&\e^2\( \ddot x +p_1(\t)\dot x\)+p_0(\t) x = f \qquad \mbox{in}\,\, (0, \ell)\\
&x(0)=x(\ell)\,\,\, \dot x(0)=\dot x(\ell)
\end{aligned}
\right.
\end{equation}
\begin{lemma}\label{lemmapbaux}
Let $$\Lambda_{p_0(\t)}=\(\int_0^\ell \sqrt{p_0(t)}\, dt\)^2.$$ There is a small number $\e_0=\e_0(p_0, \ell)>0$ such that if $\e\in (0, \e_0)$ satisfies the gap condition 
\beq\label{gap}
|4\pi^2 m^2\e^2 -\Lambda_{p_0}|\ge \tilde c_0 \e \qquad {\rm \mbox{for any}}\, m \in \mathbb N \cup \{0\}
\eeq
whit $\tilde c_0$ is small enough, then there exists a constant $C>0$ such that problem \eqref{pbaux} has a unique solution which satisfies 
\beq\label{1.1}
\e^2\|\ddot x\|_\infty +\e \|\dot x\|_\infty +\|x\|_\infty \le \frac{C}{\e}\|f\|_\infty
\eeq
for any $f\in C^0(0, \ell)$. Moreover, if in addition $f\in C^2(0, \ell)$, the unique solution to problem \eqref{pbaux} satisfies
\beq\label{1}
\e^2\|\ddot x\|_\infty +\e \|\dot x\|_\infty +\|x\|_\infty \le C\(\|\ddot f\|_\infty+\|\dot f\|_\infty+\|f\|_\infty\).
\eeq
\end{lemma}
\begin{Proof}
Although similar results were obtained in \cite{dmp}, we sketch the proof to illustrate why condition \eqref{gap} is required.\\ We take
$$\Lambda_{p_0(\t)}=\(\int_0^\ell \sqrt{p_0(t)}\, dt\)^2; \qquad s(\t)=\frac{\pi}{\sqrt{\Lambda_{p_0}}}\int_0^\t\sqrt{p_0(t)}\, dt,\qquad y(s)=x(\t). $$ Then \eqref{pbaux} is transformed into 
\begin{equation}\label{pbauxtra}
\left\{\begin{aligned}
&\e^2\( \ddot y +q(s)\dot y\)+\nu_0  y = \tilde f(s) \qquad \mbox{in}\,\, (0, \pi)\\
&y(0)=y(\pi)\,\,\, \dot y(0)=\dot y(\pi)
\end{aligned}
\right.
\end{equation}
with $$q(s)=\frac{\dot p_0}{2p_0}+\frac{1}{\pi \sqrt{\Lambda_{p_0}}}\frac{p_1}{\sqrt{p_0}};\qquad \nu_0=\frac{\Lambda_{p_0}}{\pi^2}; \qquad \tilde f(s)=\frac{{\Lambda_{p_0}}}{\pi^2}\frac{f}{p_0}.$$
It is a standard fact that the eigenvalue problem 
\begin{equation}\label{pbauxtraeig}
\left\{\begin{aligned}
&\ddot y +q(s)\dot y+\nu y = 0 \qquad \mbox{in}\,\, (0, \pi)\\
&y(0)=y(\pi)\,\,\, \dot y(0)=\dot y(\pi)
\end{aligned}
\right.
\end{equation}
has an infinite sequence of eigenvalues $(\nu_m)_m \subset \mathbb R$ such that $$\sqrt{\nu_m}=2m+O\(\frac{1}{m^3}\) \qquad \mbox{as}\,\, m\to\infty$$ with associated eigenfunctions $y_m(s)$ that forms an orthonormal basis in $L^2(0, \pi)$.\\ Thus, if $\nu_0 \neq \e^2 \nu_m$ for all $m\ge 0$ the problem \eqref{pbaux} is solvable. In such a case the solution \eqref{pbaux} can be described as $$y(s)=\sum_{m=0}^{\infty}\frac{\tilde f_m}{\nu_0-\e^2\nu_m}y_m(s)$$ where $\tilde f_m(s):=\displaystyle \int_0^\pi \tilde f(s) y_m(s)\, ds$.\\ Since $y\in C^2(0, \pi)$ the above expression holds in $C^2(0, \pi)$. From \eqref{gap} we find that $$|\nu_0 -\nu_m \e^2|\ge \frac c 2 \e$$ if $\e$ is sufficiently small. Next we notice that, by using Cauchy-Schwarz inequality and Parseval's identity we have 
\beq
\begin{aligned}
\|y\|_\infty &\le \sum_{m=0}^{\infty} \left|\frac{\tilde f_m y_m(s)}{\nu_0 - \e^2\nu_m }\right|\le \left( \sum_{m=0}^\infty \tilde f_m^2 y_m^2 \right)^{\frac 12}\left( \sum_{m=0}^\infty \frac{1}{(\nu_0 - \nu_m \e^2)^2}\right)^{\frac 12}\le  \frac{c}{\e}\|\tilde f\|_\infty 
\end{aligned}
\eeq 
Coming back to the original variable $$\|x\|_\infty =\|y\|_\infty \le\frac{ c}{\e} \left\|\frac{{\Lambda_{p_0}}}{p_0}\right\|_\infty\|f\|_\infty\le \frac{ C}{\e}\|f\|_\infty .$$ In this way, one can also estimate the $L^\infty(0, \pi)$- norms of $\dot y$ and $\ddot y$. Therefore the result holds. For a more detailed treatment of this and estimate \eqref{1} one can see [\cite{dmp}, Lemma 8.2].
\end{Proof}
In view of system \eqref{pbaux}, it is natural to consider a perturbation of the equation in \eqref{pbaux}, namely
\begin{equation}\label{pbaux2}
\left\{
\begin{aligned}
&\e^2\( \ddot x +p_1(\t)\dot x\)+\(p_0(\t)+\e\tilde p_{0, \e}(\t)\)x= f \qquad \mbox{in}\,\, (0, \ell)\\
&x(0)=x(\ell)\qquad \dot x(0)= \dot x(\ell)
\end{aligned}
\right.
\end{equation}
where $\{\tilde p_{0, \e}(\t)\}_{\e>0}$ is a family of $C^2(0, \ell)$ functions such that 
\begin{equation}\label{palpal}
\sup_{\e>0}\|\tilde p_{0, \e}\|_{C^2(0, \ell)}<C
\end{equation}
and
\begin{equation}\label{palpalpal}
\sup_{\e >0}\e \left\|\frac{\partial \tilde p_{0, \e}}{\partial \e}\right\|_\infty <C.
\end{equation}
Then we have a constant $M>0$ and a family $\{\Lambda_\e\}_\e\subset \mathbb R$ such that $$\Lambda_{p_0(\t)+\e \tilde p_{0, \e}(\t)}=\Lambda_{p_0(\t)}+ \e \Lambda_\e$$ and 
\begin{equation}\label{4.19}
|\Lambda_\e|+\e\left|\frac{\partial \Lambda_\e}{\partial \e}\right|\le M.
\end{equation}
We observe that if there exists a small $\e>0$ such that 
\begin{equation}\label{gap1}
|4\pi^2 m^2 \e^2 -\(\Lambda_{p_0}+\e \Lambda_\e\)|\ge \tilde c_0 \e \qquad m=0, 1, 2, \ldots
\end{equation}
for some small $\tilde c_0>0$, then \eqref{gap} holds afetr $\Lambda_{p_0}$ is substituted by $\Lambda_{p_0+\e \tilde p_{0, \e}}$ and hence existence of a unique solution to \eqref{pbaux2} satisfying a priori bounds \eqref{1.1} and \eqref{1} is guaranteed.\\ Moreover \eqref{palpal} allows us to choose the constant $C>0$ in \eqref{1.1} and \eqref{1} to be independent of $\e$.
\begin{remark}\label{4.3}
\begin{itemize}
\item[(i)]{\rm First we deduce a sufficient condition of $\e>0$ for which inequality \eqref{gap1} holds. Notice that \eqref{gap1} means that if} $$4\pi^2m^2\e^2 \ge \Lambda_{p_0}+\e\Lambda_\e \qquad (\mbox{or}\,\, 4\pi^2 m^2\e^2 \le \Lambda_{p_0}+\e \Lambda_\e)$$ {\rm then it should be} $$4\pi^2m^2\e^2- \tilde c_0 \e \le \Lambda_{p_0}+\e\Lambda_\e \qquad (\mbox{or}\,\, 4\pi^2m^2\e^2+\tilde c_0\e \le \Lambda_{p_0}+\e\Lambda_\e)$$
{\rm for a sufficiently small $c>0$ and for every $m\in\mathbb N \cup \{0\}$. Given any small number $\e>0$, let us write }
\begin{equation}\label{bla}
4\pi^2\e^2 =\frac{\Lambda_{p_0}+\e \Lambda_\e}{(m_0+a_0)^2}
\end{equation}
{\rm with some $m_0\in\mathbb N$ large and $a_0\in [0, 1)$. Assume $a_0\neq 0$. Then the least $m\in\mathbb N$ satisfying $4\pi^2m^2\e^2\ge \Lambda_{p_0}+\e\Lambda_\e$ is $m=m_0+1$. Besides, for $m\ge m_0+1$ we have }
\begin{equation*}
\begin{aligned}
4\pi^2m^2\e^2-\tilde c_0\e &\ge \(\Lambda_{p_0}+\e\Lambda_\e\)\left[\left(\frac{m_0+1}{m_0+a_0}\right)^2-\tilde c_0 \frac{\(\Lambda_{p_0}+\e\Lambda_\e\)^{-\frac 12}}{2\pi(m_0+a_0)}\right]\\
&\ge \(\Lambda_{p_0}+\e\Lambda_\e\)\left[1+\frac{2(1-a_0)}{m_0}-\frac{\tilde c_0}{2\pi \sqrt{\Lambda_{p_0}}}\frac{1}{m_0}+o\(\frac{1}{m_0^2}\)\right]\\
&\ge \Lambda_{p_0}+\e\Lambda_\e
\end{aligned}
\end{equation*}
{\rm provided $a_0 \le 1-\frac{\tilde c_0}{2\pi\sqrt{\Lambda_{p_0}}}$ choosing $\tilde c_0<2\pi\sqrt{\Lambda_{p_0}}$.}
\item[(ii)] {\rm Let us show the existence of a sequence of small positive numbers $\e>0$ converging to zero satisfying \eqref{gap1} provided \eqref{palpalpal} holds.\\ Indeed it is easy to see that the equation \eqref{bla} has a unique pair $(m_0, a_0)$ for any $\e\in (0, \e_1)$ where $\e_1>0$ is determined by $\Lambda_{p_0}$ and $M$ in \eqref{4.19}.}
\end{itemize}

\end{remark}
We come back to the original problem.\\
Let us introduce the linear operator
$$L_0(e_0):=\e^2 \(a_0(\e s) \ddot e_0 +a_1(\e s) \dot e_0\)+a_2(\e s) e_0.$$
The following result holds.
\begin{lemma}\label{invL}
 We have a positive number $\Lambda_{\hat \mu_0}$ and a number $\{\Lambda_{\hat \mu_0, \e}\}_\e$ such that if
\begin{equation}\label{risonanza}
|4\pi^2 m^2\e^2 -(\Lambda_{\hat \mu_0}+\e \Lambda_{\hat \mu_0, \e})|\ge \tilde c_0 \e \quad m=1, 2\ldots \eeq
 for some positive and sufficiently small constant $\tilde c_0$, then for any $f\in C^0_\ell(\mathbb R)\cap L^\infty (\mathbb R)$, there exists a unique $e_0\in C^2_{\ell}(\mathbb R)$ solution of $L_0(e_0)=f$. Moreover, there exists $C>0$ such that $$\|e_0\|_\e=\e^2\|\ddot e_0\|_\infty+\e\|\dot e_0\|_\infty+\|e_0\|\leq \frac{C}{\e}\|f\|_\infty.$$ Finally, if $f\in C^2_\ell (\mathbb R)$, then $$\|e_0\|_\e=\e^2\|\ddot e_0\|_\infty+\e\|\dot e_0\|_\infty+\|e_0\|\leq C\[\|f\|_\infty+\|\dot f\|_\infty+\|\ddot f\|_\infty\].$$

\end{lemma}
\begin{Proof}
The equation $\e^2\(a_0(\e s)\ddot e_0+a_1(\e s)\dot e_0\)+a_2(\e s) e_0 =f$ can be written as 
$$\e^2 \(\ddot e_0 + p_1(\t) \dot e_0\)+(p_0(\t)+\e \tilde p_{0, \e}) e_0 = g$$ with $$p_1(\t)=\frac{a_1(\e s)}{a_0(\e s)};\qquad p_0(\t)=\frac{\Lambda_1}{\hat\mu_0^2(\t)};\qquad \tilde p_{0, \e}(\t)=\frac{a_2^\e(\t)\hat\mu_0^2(\t)-a_0^\e\Lambda_1}{\hat\mu_0^4(\t)}+\e q_\e(\t);\qquad g=\frac{f(\t)}{a_0(\t)}.$$ It is clear that $p_0(\t)>0$ is a $C^2(0, \ell)$ function and $\tilde p_{0, \e}\in C^2(0, \ell)$ function and \eqref{palpal} and \eqref{palpalpal} hold. Then we let $$\Lambda_{\hat\mu_0}=\left(\int_0^\ell \sqrt{\frac{\Lambda_1}{\hat\mu_0^2}}\, dt\right)^2$$ and hence there exist numbers $\Lambda_{\hat\mu_0, \e}$ such that $$\Lambda_{p_0+\e\tilde p_{0, \e}}=\Lambda_{\hat \mu_0}+\e \Lambda_{\hat\mu_0, \e}$$ and the result comes from the above discussions. 
\end{Proof}
\begin{Proof}[Proof of Theorem \ref{principale}]
By   Lemma \ref{invL} it follows that there exists a sequence of small $\e=\e_m>0$ converging to zero as $m\to+\infty$  such that the operator $L_0(e_0)$ is invertible with bounds for $L_0(e_0)=h$ given by $$\| e_0\|_\e\leq C\e^{-1}\|h\|_\infty ,$$ for some positive constant $C.$ Finally, by Contraction Mapping Argument using the properties of the right-hand side of \eqref{sistemaode}, it follows that, the problem \eqref{sistemaode} has a unique solution with
$$\|  e_0\|_\e< c \e^{(\sigma+2)(1-a)}
$$
and that concludes the proof.
\end{Proof}

\section{The linear theory}
\label{4}
In this section we give the proof of Proposition \ref{inv}.
 We need a couple of preliminary results.

\begin{lemma}\label{lin1}
Assume $\xi\not\in\{0, \pm \sqrt{\Lambda_1}\}$. Then given $h\in L^\infty(\mathbb R^2)$, there exists a unique bounded solution of
\beq\label{pblin1}
(   \hat {\mathcal L}-\xi^2)\psi=h \ \hbox{in $\mathbb R^2$.}
\eeq
Moreover
\beq\label{stimapriori}
\|\psi\|_\infty \leq C_{\xi}\|h\|_\infty
\eeq
for some constant $C_\xi>0 $ only depending on $\xi$.
\end{lemma}
\begin{Proof}We argue as in Lemma 3.1 of \cite{dmp}. \end{Proof}\\\\
  \begin{lemma}\label{sollim}
Let $\phi $ a bounded solution of $\hat  {\mathcal L}(\phi)+\partial^2_{ss}\phi=0$ in $\mathbb R^2$. Then $\phi(s, t)$ is a linear combination of the functions $Z_1(t)$, $Z_0(t)\cos(\sqrt{\Lambda_1}s)$, $Z_0(t)\sin(\sqrt{\Lambda_1}s)$.
\end{lemma}
\begin{Proof}
 We argue as in Lemma 7.1 of \cite{delpinotoda}.
 \end{Proof}
 \\\\
 \noindent
 {\bf Proof of Proposition \ref{inv}}
The proof  will be carried out in three steps.\\
\\
{\em Step 1: A priori bound (special case)}     Let us assume for the moment that in problem \eqref{gluing21} the function $c_0$ is identically zero.\\ We will prove that there exits $C>0$ so that for any $h$ with $\|h\|_{**}<+\infty$ and any $\phi$ solution of problem
\begin{equation}\label{gluing21primo}
\left\{
\begin{array}{lr}
\mathcal L(\phi)= h\qquad\qquad \qquad\mbox{in}\,\,\, \mathcal C_{2\d}\\\\
\partial_\nu \phi=0\qquad\qquad\quad\qquad \mbox{on}\,\,\, \partial\mathcal C_{2\d}\cap \{t=0\}\\\\
\phi\(s+\frac{\ell}{\e}, t\)=\phi(s, t)\\\\
\displaystyle\int_{-\frac{2\d}{\mu}}^0 \phi(s, t)Z_0(t)\, dt=0\qquad \forall\,\, s\in \mathbb R^+.
\end{array}
\right.
\end{equation}

 with $\|\phi\|_{*}<+\infty$ we have $$\|\phi\|_{*}\leq C \|h\|_{**}.$$
By contradiction we assume that there exist sequences $\l_n\rightarrow 0$, $(h_n)_n$ and $(\phi_n)_n$ solutions of
\eqref{gluing21primo} where $$\d_n=\e_n^a\ \hbox{for some $a\in(0,1)$ and}\  \mu_n(\e_n s)=\e_n \hat\mu(\e_n s)$$ such that $$\|\phi_n\|_{*}=1\qquad\qquad \|h_n\|_{**}\rightarrow 0.$$
 To achieve a contradiction we will first show that

\begin{equation}\label{normainftyphi}
\|\phi_n\|_{\infty}\rightarrow 0.
\end{equation}
If this was not the case then we may assume that there is a positive number $c$ for which $\|\phi_n\|_{\infty,  \mathcal C_{2\d_n}}>c$. Since we also know that $$|\phi_n(s, t)|\leq \frac{c}{(1+|t|)^\sigma},$$ we conclude that for some $A>0$ $$\|\phi_n\|_{L^\infty(|t|\leq A)}\geq c.$$ Let us fix an $s_n$ such that $$\|\phi_n(s_n, \cdot)\|_{L^\infty(|t|\leq A)}\geq \frac{c}{2}.$$ By elliptic estimates, compactness of Sobolev embeddings and the fact that the coefficients of $\tilde{\mathcal A}(\phi_n)$ tends to zero as $\l_n\rightarrow 0$, we see that we may assume that the sequence of functions $\tilde\phi_n(s,t):=\phi_n(s+s_n, t)$ converges uniformly over compact subsets of $\mathbb R^2$, to a nontrivial, bounded solution of $$\hat\mu_0^\infty\partial_{ss}^2\tilde{\phi}+\partial_{tt}^2\tilde{\phi}
+e^w\tilde{\phi}=0\qquad\mbox{in}\,\,\, \mathbb R^2$$
where $\hat\mu_0^\infty$ is a positive constant, which with no loss of generality via scaling, we may assume equal to one. By virtue of Lemma \ref{sollim} then $\tilde \phi$ is a linear combination of $Z_0$ and $Z_1$. Moreover    { by the decay behavior and the orthogonality conditions assumed, which pass to the limit thanks to the Dominated Convergence, we find then that $\tilde{\phi}\equiv 0$. }This is a contradiction that shows the validity of \eqref{normainftyphi}.\\\\ Let us conclude now the result of Step 1. \\ Since $\|\phi_n\|_{*}=1$, there exists $(s_n, t_n)$ with $r_n:=|t_n|\rightarrow +\infty$ such that $$r_n^\sigma |\phi_n(s_n, t_n)|+r_n^{\sigma+1}|D\phi_n(s_n, t_n)|\geq c>0.$$ Let us consider now the scaled function $$\tilde{\phi}_n(z_0, z)=r_n^\sigma \phi_n (s_n+r_n z_0, r_n z)$$ defined on $\bar{\mathcal D}$ given by
$$\bar{\mathcal D}:=\left\{(z_0, z)\,\, :\,\,-r_n^{-1}s_n \leq z_0 \leq r_n^{-1}\left(\frac{\ell}{\e_n}-s_n\);\quad  -\frac{2\d_n r_n^{-1}}{\mu_n(\e_n s)}<z<0\right\}.$$
Then we have $$|\tilde{\phi}_n(z_0, z)|+|z||D\tilde{\phi}_n(z_0, z)|\leq |z|^{-\sigma}\qquad \mbox{in}\,\, \bar{\mathcal D}$$ and for some $z_n$ with $|z_n|=1$ $$|\tilde{\phi}_n(0, z_n)|+|D\tilde{\phi}_n(0, z_n)|\geq c>0.$$
 Moreover $\tilde{\phi}_n$ satisfies $$\tilde{\mu}_{0, n}^2 \partial_{z_0z_0}^2\tilde{\phi}_n+\partial_{zz}^2\tilde{\phi}_n+o(1) \tilde{C}(\tilde{\phi}_n)=\tilde{h}_n\qquad \mbox{in}\,\, \bar{\mathcal D}$$
 where $\tilde{h}_n(z_0, z)=r_n^{\sigma+2}h_n(s_n+r_n z_0, r_n z)$, $\tilde{\mu}_{0, n}=\hat\mu_0^n(s_n+r_nz_0)$ and $\tilde{C}(\tilde{\phi}_n)$ is bounded. \\
Since $$\|\partial_s\hat\mu_{0, n}^2\|_{\infty, \bar
D}=O(\e_n);\qquad \left\|\partial_s\(\frac{r_n^{-1}s_n}{\mu_n(\e_n(s_n+r_n
z_0))}\)\right\|_{\infty, \bar D}=O(\e_n^a);$$ $$\qquad
\left\|\partial_{ss}^2\(\frac{r_n^{-1}s_n}{\mu_n(\e_n(s_n+r_n
z_0))}\)\right\|_{\infty, \bar D}=O(\e_n^{1+a})$$ then, we may assume that

 $$\tilde{\mu}_{0, n}^2 \rightarrow \hat\mu^*>0$$ and that the function $\tilde{\phi}_n$ converges uniformly, in $C^1-$ sense over compact subsets of $\mathcal D_*$, to $\tilde\phi$ which satisfies
 \begin{equation}\label{eqtildephi}
 \hat\mu^* \partial_{z_0z_0}\tilde\phi+\partial_{zz}\tilde\phi=0\qquad \mbox{in}\,\, \mathcal D_*
 \end{equation}
 where $$\mathcal D_*:=(0, \infty)\times (-\infty, 0)$$ and $\tilde\phi$ satisfies
 \begin{equation}\label{tildephi}
 |\tilde\phi(z_0, z)|+|z||D\tilde\phi(z_0, z)|\leq |z|^{-\sigma}
 \end{equation}
 with the boundary condition. With no loss of generality, we may assume that $\hat\mu^*=1$.\\
Hence $\tilde\phi$ is weakly harmonic in $\mathcal D_*$ and hence $\tilde \phi \equiv const$. Moreover since it satisfies \eqref{tildephi}, it follows that $\tilde\phi\equiv 0$.\\ This is a contradiction.\\\\

  {\em  Step 2: A priori bound (general case) } We claim that the a priori estimate obtained in Step 1 is valid for the full problem \eqref{gluing21}. We conclude from Step 1 that 
\begin{equation}\label{A1}
\|\phi\|_* \le c \left[\|h\|_{**}+ \|c_0 Z_0\|_{**}\right]\le C \left[\|h\|_{**}+ \|c_0\|_\infty\right]
\end{equation}
for any $h$ with $\|h\|_{**}<\infty$ and solution $\phi$ of problem \eqref{gluing21}. To conclude we have to find a bound for the coefficient $c_0(s)$.\\
   Testing the equation in \eqref{gluing21} with $Z_0$ and integrating with respect to
   $dt$, we get
 \begin{equation}\label{testZ1}
 c_0( s)\int_{-\frac{2\d}{\mu}}^0Z_0^2\, dt = \int_{-\frac{2\d}{\mu}}^0\mathcal L (\phi) Z_1\, dt -\int_{-\frac{2\d}{\mu}}^0 h Z_1\,
 dt
 \end{equation}

Since $Z_0$ decays exponentially
 $$\int_{-\frac{2\d}{\mu}}^0 Z_0^2\, dt =\frac 1 2+
 O(e^{-\sqrt{\Lambda_1 }\frac{\d}{\mu}})$$
 Hence from \eqref{testZ1} it follows that 
 \begin{equation}\label{testZ1primo}
\begin{aligned}
 c_0(s)\(\frac 12 + O\(e^{-\sqrt{\Lambda_1}\frac{1}{\e}}\)\)&=\int_{-\frac{2\d}{\mu}}^0 \hat\mu_0^2 \partial_{ss}^2 \phi Z_0\, dt +\int_{-\frac{2\d}{\mu}}^0 \(\partial^2_{tt}\phi+e^w\phi\) Z_0\, dt+\int_{-\frac{2\d}{\mu}}^0 \tilde{\mathcal A}(\phi) Z_0\, dt \\
&-\int_{-\frac{2\d}{\mu}}^0 h Z_0\,
 dt.
\end{aligned}
 \end{equation}

 It is easy to see that
 \beq\label{stimah}
 \left|\int_{-\frac{2\d}{\mu}}^0 h Z_0\, dt\right|\leq
 \|h\|_{**}\int_{-\frac{2\d}{\mu}}^{0}\frac{1}{(1+|t|)^{\sigma+2}}\,
 dt \leq C\|h\|_{**}.
 \eeq

 Now by using the boundary condition, the orthogonality condition and the radial symmetry of $Z_0$
 we get
 \begin{equation}\label{linZj}
\left|\int_{-\frac{2\d}{\mu}}^{0}(\partial_{tt}^2 \phi + e^w \phi)
Z_0\, dt\right|\leq  O(e^{-\frac{c}{\e}})\|\phi\|_*.
 \end{equation}

 Now,  since
  $$\int_{-\frac{2\d}{\mu}}^0\phi Z_0\, dt =0$$ if we make twice the $s-$ derivative and we make some computations, we immediately get that
  $$\left|\int_{-\frac{2\d}{\mu}}^0\partial_{ss}^2\phi Z_j\, dt\right|\le O\(e^{-\frac c \e}\)\|\phi\|_*.$$

  Moreover we get
\begin{eqnarray*}
\left|\int_{-\frac{2\d}{\mu}}^0 \tilde{\mathcal A}(\phi)Z_0\, dt\right| &=&\left|\int_{-\frac{2\d}{\mu}}^0 b_0(s, t) \partial_{ss}^2\phi Z_0\, dt+\int_{-\frac{2\d}{\mu}}^0 b_1(s, t)\partial_{tt}^2 \phi Z_0\, dt\right. \\
&&\left.+\int_{-\frac{2\d}{\mu}}^0 b_2(s, t)\partial^2_{st}\phi Z_0\, dt+\int_{-\frac{2\d}{\mu}}^0 b_3(s, t) \partial_t\phi Z_0, dt +\int_{-\frac{2\d}{\mu}}^0 b_4(s, t) \partial_s\phi Z_0\, dt\right|
\end{eqnarray*}
where $b_j$ are defined in \eqref{tildeA}.\\
Now reasoning as before
\begin{eqnarray*}
\left|\int_{-\frac{2\d}{\mu}}^0   b_0(s, t) \partial_{ss}^2\phi Z_0 \, dt \right|&\leq & \e \int_{-\frac{2\d}{\mu}}^0 |\partial_{ss}^2\phi Z_0|\, dt \leq  \e O\(e^{-\frac c \e}\) \|\phi\|_*
\end{eqnarray*}

\begin{equation*}
\left| \int_{-\frac{2\d}{\mu}}^0 b_1(s, t) \partial_{tt}^2\phi Z_0\, dt \right|\leq  O(e^{-\frac c \e})\|\phi\|_*+C \e^2\|\phi\|_*\le C \e^2\|\phi\|_*
\end{equation*}

\begin{equation*}
\left|\int_{-\frac{2\d}{\mu}}^0 b_2(s, t) \partial_{st}^2\phi Z_0 \, dt \right|\leq \e \|\phi\|_*;\qquad 
\left| \int_{-\frac{2\d}{\mu}}^0 b_3(s, t) \partial_{t}\phi Z_0 \, dt \right|\leq  \e \|\phi\|_*; \qquad
\end{equation*}
\begin{equation*}
\left| \int_{-\frac{2\d}{\mu}}^0 b_4(s, t) \partial_{s}\phi Z_0 \, dt \right|\leq   \e^2 \|\phi\|_*
\end{equation*}

hence by \eqref{testZ1primo}
\begin{equation}\label{A}
\|c_0\|_\infty\leq C\|h\|_{**}+c\e\|\phi\|_*.
\end{equation}

Combining \eqref{A} with \eqref{A1} the result follows.\\

 {\em Step 3: (Existence part)} We establish now the existence of a solution $\phi$ for problem \eqref{gluing21}.\\ We consider the case in which $h(s, t)$ is a $T$-periodic function in $s$, for an arbitrarily and large but fixed $T$. We then look for a weak solution $\phi$ to \eqref{gluing21} in $H_T$ defined as the subspace of functions $\psi$ which are in $H^1(B)$ for any $B$ bounded subset of $\mathcal C_{2\d}$, which are $T$-periodic in $s$, such that $\partial_\nu\phi=0$ on $\partial \mathcal C_{2\d}\cap \{t=0\}$ and so that $$\int_{-\frac{2\d}{\mu}}^0 \psi Z_0\, dt=0\qquad \forall\,\,\ \psi\in H^1(B).$$
 Let $\mathcal D_T:=\{ t\in\left[-\frac{2\d}{\mu}, 0\right]\,\,:\,\, s\in (0, T)\}$ and the bilinear form in $H_T$: $$\mathcal {B}(\phi, \psi):=\int_{{\mathcal D}_T} \psi \mathcal{L}(\phi)\, dt\qquad \forall\,\, \psi\in H_T.$$ Then problem \eqref{gluing21} gets weakly formulated as that of finding $\phi\in H_T$ such that $$\mathcal{B}(\phi, \psi)=\int_{\mathcal D_T} h\phi\, dt\qquad \forall\,\, \psi\in H_T.$$ If $h$ is smooth, elliptic regularity yields that a weak solution is a classical one.\\ The weak formulation can be readily put into the form $$\phi+\mathcal{K}(\phi)=    {h}\qquad \mbox{in}\,\, H_T$$ where $     h$ is a linear operator of $h$ and $\mathcal K$ is compact.\\ The a priori estimate of Step 2 yields that for $h=0$ only the trivial solution is present. Fredholm alternative thus applies yielding that problem \eqref{gluing21} is thus solvable in the periodic setting. This is enough for our purpose. However we remark that if we approximate a general $h$ by periodic functions of increasing period and we use uniform estimate we obtain in the limit a solution to the problem.

\bigskip
\noindent
{\bf Acknowledgements.} {
The first  author has been supported by grants Fondecyt  1150066, Fondo Basal CMM and by Millenium Nucleus CAPDE NC130017. The second and the third authors have been partially supported by MIUR-PRIN project-201274FYK7-005 and by    the group GNAMPA of Istituto Nazionale di Alta Matematica (INdAM).}

\end{document}